\documentclass[]{article}
\usepackage{graphicx,subfigure}
\usepackage{amsfonts,amsbsy,amscd,amsgen,dsfont,amsmath,amssymb,mathrsfs,amsthm}
\usepackage{bm,multirow,esvect,cancel,xcolor,comment}
\usepackage{enumitem} 
\usepackage[colorlinks=true,allcolors=blue]{hyperref}
\usepackage{authblk}
\usepackage{algorithm}
\usepackage{algpseudocode} 
\newcounter{algsubstate}
\renewcommand{\thealgsubstate}{\alph{algsubstate}}
\newenvironment{algsubstates}
{\setcounter{algsubstate}{0}%
	\renewcommand{\State}{%
		\stepcounter{algsubstate}%
		\Statex {\footnotesize\thealgsubstate:}\space}}
{}
\usepackage[utf8]{inputenc}
\usepackage[T1]{fontenc}

\newcommand{\id}{\mathrm d}

\newcommand{\etab}{\pmb{\eta}}
\newcommand{\bxi}{\pmb{\xi}}

\renewcommand{\tilde}{\widetilde}

\DeclareMathAlphabet\mathbfcal{OMS}{cmsy}{b}{n}

\newcommand{\R}{\mathbb{R}}

\renewcommand{\a}{\alpha} 
\graphicspath{{./figures-v2/}}

\newcommand{\sD}{\mathcal D}
\newcommand{\sS}{\mathcal S}
\newcommand{\grad}{\nabla} 

\renewcommand{\d}{\partial}
\newcommand{\g}{\gamma} 
\renewcommand{\b}{\beta}

\newcommand{\cb}{} 

\newtheorem{theorem}{Theorem}

\begin{document}

\title{Data-driven prediction of multistable systems from sparse measurements}
\author[]{Bryan Chu} 
\author[]{Mohammad Farazmand\thanks{Corresponding author's email address: farazmand@ncsu.edu}}
\affil{Department of Mathematics, North Carolina State University,
	2311 Stinson Drive, Raleigh, NC 27695-8205, USA}
\date{}

\maketitle

\begin{abstract}
    We develop a data-driven method, based on semi-supervised classification, to predict the asymptotic state of multistable systems when only sparse spatial measurements of the system are feasible. 
Our method predicts the asymptotic behavior of an observed state by quantifying its proximity to the states in a precomputed library of data.
To quantify this proximity, we introduce a sparsity-promoting metric-learning (SPML) optimization, which learns a metric directly from the precomputed data. The optimization problem is designed so that
the resulting optimal metric satisfies two important properties: (i) It is compatible with the precomputed library, and (ii) It is computable from sparse measurements. We prove that the proposed SPML optimization is convex, its minimizer is non-degenerate, and it is equivariant with respect to scaling of the constraints.
We demonstrate the application of this method on {\cb two multistable systems:} a reaction-diffusion equation, arising in pattern formation, which has four asymptotically stable steady states {\cb and a FitzHugh–Nagumo model with two asymptotically stable steady states}. 
Classifications of the multistable reaction-diffusion equation based on SPML predict the asymptotic behavior of initial conditions based on two-point measurements with $95\%$ accuracy
when moderate number of labeled data are used. {\cb For the FitzHugh-Nagumo, SPML predicts the asymptotic behavior of initial conditions from one-point measurements with $90\%$ accuracy.}
The learned optimal metric also determines where the measurements need to be made to ensure accurate predictions.
\end{abstract}


\section{Introduction}\label{sec:Introduction}
Predicting spatiotemporal dynamics, described by partial differential equations (PDEs), requires the initial state of the system to be
measured on a dense spatial grid. In practice, such a dense set of measurements is often unavailable. 
This problem is specially conspicuous in fluid turbulence~\cite{Farazmande1701533,Westerweel2013}, meteorology~\cite{Tomlinson2011} and oceanography~\cite{Lumpkin2017} where the dense direct observations are out of reach.
When only sparse spatial measurements of the system are available, the PDE model cannot be numerically integrated to predict the evolution
of the system. A common remedy is to use data assimilation techniques, such as 4DVAR and ensemble Kalman filtering, to estimate the state of 
the system from sparse observational data~\cite{bannister2017,GHIL1991,Stuart2015}. 
These data assimilation techniques are computationally expensive, and therefore not ideal for making predictions in real-time. 

\begin{figure}
	\centering
	\includegraphics[width=.9\textwidth]{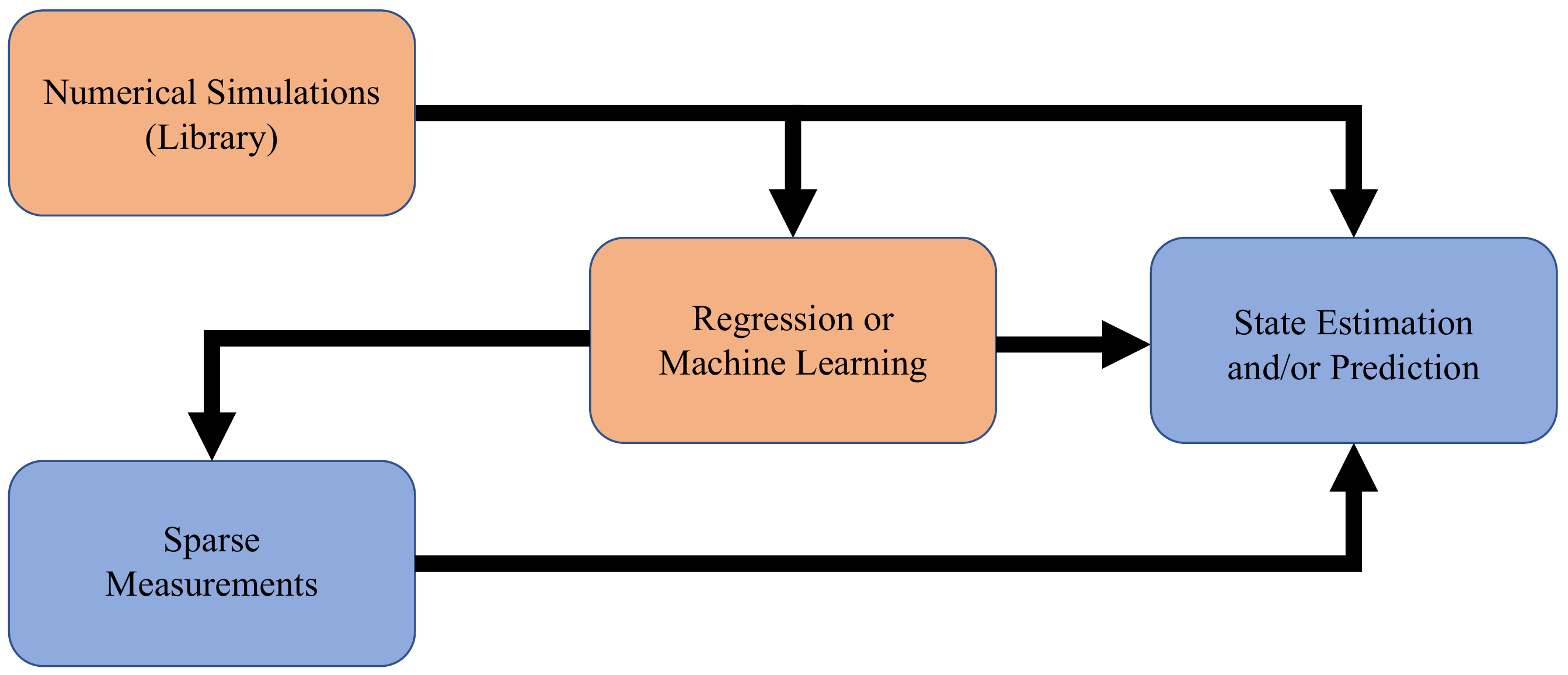}
	\caption{Summary of the algorithm used for data-driven state estimation and/or prediction. A large library of system states and dynamics is generated offline by numerical simulations. 
		This library is then fed into an algorithm (regression or machine learning block) to inform the sparse measurements. Finally, the library and sparse measurements are blended to estimate the 
		system state or to make predictions. The blocks in orange contain the offline computations, whereas the blocks in blue are performed in real time. }
	\label{fig:schem_algorithm}
\end{figure}

An alternative data-driven approach is to estimate the system state from a library of precomputed states. These methods shift the bulk of the computational cost
offline, speeding up the real-time predictions or state estimation (see section~\ref{sub:related_work} for a review).
As summarized in figure~\ref{fig:schem_algorithm}, this approach seeks to interpolate a large library of precomputed states 
to find a state estimation that best agrees with the sparse observational data.

As a critical component, these methods require a measure of proximity among the precomputed states, and between precomputed states and the observational data.
The proximity between states is measured using a distance metric which is induced by a functional norm. This naturally raises two important questions:
\begin{enumerate}
\item What is the best metric to use?
\item Is the metric computable from sparse measurements?
\end{enumerate}
Here, we develop a method that addresses both these questions in the context of multistable dynamical systems, which have multiple co-existing attractors.
Based on the initial condition, the system evolves towards one of these attractors. 
Our goal is to predict the asymptotic state of the system given sparse measurements of its initial condition.
To this end, we create a library of initial states from offline numerical simulations, and label the states according to the attractor to which they converge. 
Given the sparse measurements of a new system state, we predict its asymptotic behavior by classifying it using the information from the precomputed library. 
This classification requires a measure of similarity between the labeled states and the observational data. 
This similarity is quantified using a metric which is induced by a norm.

To find an optimal norm that is quantifiable from sparse measurements, we formulate and solve a sparsity-promoting metric-learning (SPML) optimization problem. 
The optimization problem is designed so that similarly labeled data are relatively close when their distance is measured in the optimal norm. 
At the same time, the optimization problem requires that differently labeled data have a relatively large distance.
Most importantly, a sparsity-promoting penalization is added to the cost function to ensure that the optimal norm is quantifiable from
sparse measurements.

\subsection{Related Work}%
\label{sub:related_work}

In this section, we review the previous studies that use a library of precomputed data for parameter or state estimation 
from sparse observational data.

Bright et al.~\cite{Bright} developed a method for reconstructing the pressure field in flow past a cylinder.
To this end, they solve a convex $L^1$ optimization problem to 
find a sparse representation of the flow from a precomputed library of POD modes. 
Kramer et al.~\cite{Kramer} take a similar approach but instead of POD modes, they use a library obtained from 
Dynamic Mode Decomposition (DMD). They use this DMD library together with sparse sensing in order to classify different flow regimes (e.g., steady, periodic, or chaotic)
in a thermofluid system.
Callaham et al. \cite{Callaham2019} report smaller state estimation errors 
by using a library of flow snapshots as opposed to using a library of model modes, such as POD or DMD modes.

In the above-mentioned studies, the sparsity-promoting optimization is used to minimize the number of 
library elements used to estimate the system state, not to optimize the sensor locations.
In contrast, Manohar et al.~\cite{Manohar} focus on the optimal location of the sensors. They compare several methods for
optimal sensor placement, including randomized sensing and QR pivoting. They show that flow reconstruction can be improved
by judicious choice of sensor locations.
Brunton et al.~\cite{BruntonBing2016} combine linear discriminant analysis (LDA) and principle component analysis (PCA) to determine the optimal sensor placement to be used in a classification 
algorithm. Their method is largely successful and computationally inexpensive; however, it requires the number of sensors to be predetermined by the user. In contrast, our method determines the optimal number of sensors as well as their optimal locations at the same time. Furthermore, our approach allows for the use of general classes of metrics.

Above references invariably use ad hoc metrics to measure the similarity between states; these are usually the $\ell^2$ norm applied to a discretization of the continuous state.
In contrast, here we seek to learn an optimal metric given the available data.
The importance of the metric has long been recognized in data science and machine learning communities~\cite{Coifman2005,Colonius2002,Kevrekidis2018,Kevrekidis2006,Weinberger,Xing2003}. 
In fact, our method draws on the work of Xing et al.~\cite{Xing2003} who proposed a metric learning algorithm for semi-supervised clustering of data in finite dimensions.
We extend this approach to infinite-dimensional function space and introduce a sparsity-promoting $L^1$ regularization to the optimization in order to enable the metric evaluation from sparse sampling of the functions. Furthermore, our formulation introduces an additional $L^2$ regularization to the optimization problem. We show that, without this penalization, the method of Xing et al.~\cite{Xing2003} generically returns rank-1 metrics which are too restrictive.

\subsection{Organization of this paper}
In section~\ref{sec:Sparse}, we describe the set-up of the prediction problem for multistable systems.
In section~\ref{sec:opt}, we introduce our sparsity-promoting metric-learning optimization problem.
Section~\ref{sec:Numerical Results} contains the application of our method on two examples: a reaction-diffusion PDE {\cb and the FitzHugh--Nagumo equation.}
We present the concluding remarks in section~\ref{sec:Conclusion}.

\section{Problem set-up: multistable systems}
\label{sec:Sparse}
%

We consider multistable dynamical systems, described by the general partial differential equation,
\begin{equation}
    \frac{\d u }{\d t}   = \mathcal F(u),\quad u(x,t_0)=u_0(x),
    \label{eq:pde}
\end{equation}
where $\mathcal F$ is a differential operator and $u: \Omega \times \mathbb R^+\to \mathbb R^m$ denotes the state of the system
over the spatial domain $\Omega\subset\mathbb R^n$. The initial state of the system at time $t=t_0$ is given by $u_0$.

\begin{figure}
	\centering
	\includegraphics[width=.9\textwidth]{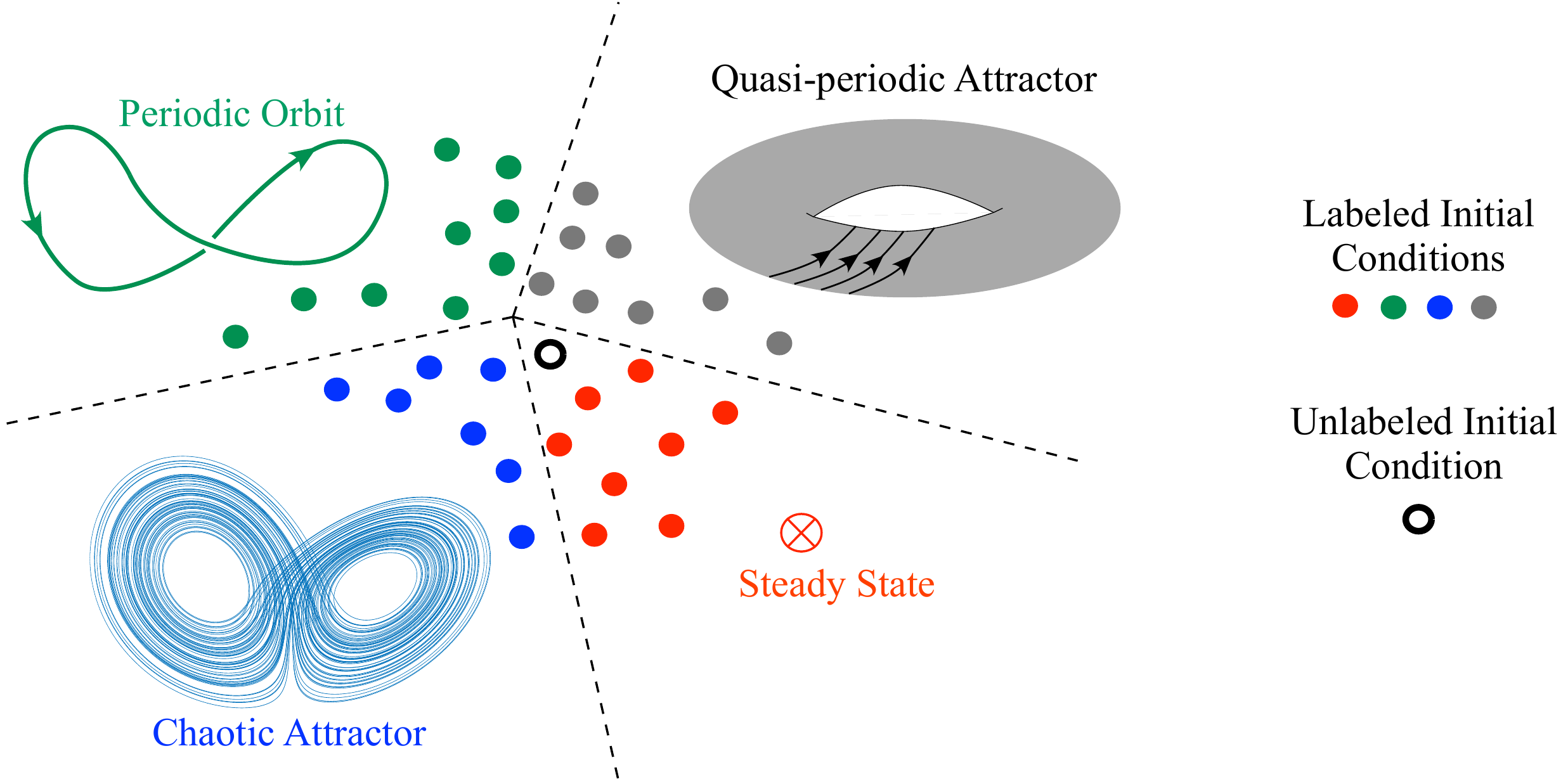}
	\caption{Schematic description of the classification problem for multistable systems. 
		Each attractor has its own basin of attraction whose boundary is marked by dashed lines. The filled circles represent the 
		initial conditions, labeled according to their asymptotic behavior. The empty black circle represents the new initial condition to be classified based on the labeled data.}
	\label{fig:schem_multistable}
\end{figure}

We assume that the system is multistable, i.e., it has $k$ asymptotically stable attracting sets $\mathcal A_1, \ldots, \mathcal A_k$.
The attractors may be steady states (or equilibria), periodic orbits, quasi-periodic orbits, or chaotic attractors (see figure~\ref{fig:schem_multistable}, for an illustration). 
Multistable systems are ubiquitous in nature; examples include tipping points in climatology~\cite{farazmand2019c,Lenton2008,Scheffer2012}, and
pattern formation in biology~\cite{gierer1972}, chemistry~\cite{Cross1993,maini1997}, and physics~\cite{gollub1999}. 
While our method applies to the general case depicted in figure~\ref{fig:schem_multistable} (coexisting attractors of different types), in practice the attractors usually comprise 
steady states and/or periodic orbits.

Given the initial state $u_0(x)$ of the system, we seek to predict its asymptotic behavior, i.e., the attractor $\mathcal A_\beta$ to which it converges as $t\to \infty$.
Of course, if the initial condition is available on a dense spatial grid over $\Omega$, PDE~\eqref{eq:pde} can be numerically integrated to determine the asymptotic behavior.
Here, however, we consider the more practical case where the initial state is only measured on a set of sparse discrete points $x_j\in \Omega$, $j=1,2,\ldots, n_s$, 
where $n_s$ denotes the number of sensors. We denote the set of sparse sensor locations by $\mathcal G_s = \{x_1,x_2,\ldots,x_{n_s}\}$.

When the number of sensors $n_s$ is small, the PDE cannot be numerically integrated. Instead, we generate an offline library of states by taking $N$ initial conditions
$u^{(i)}_0(x)$, $i=1,2,\ldots,N$, defined on a dense spatial grid. These initial conditions are in turn evolved under PDE~\eqref{eq:pde} to determine their 
asymptotic behavior. We assign a label $\ell(i)$ to the $i$-th initial condition, so that $\ell(i) = \beta\in \{1,2,\ldots,k\}$ if the initial condition $u_0^{(i)}$ converges to attractor $\mathcal A_\beta$
as $t\to \infty$. This constitutes our precomputed library of labeled data.

Now given sparse measurements of a new initial condition $u_0(x_j)$, where $x_j\in\mathcal G_s$, we would like to predict its asymptotic behavior by leveraging the 
precomputed library of labeled data. This set-up reduces our prediction problem to a semi-supervised classification problem where the labeled data are numerically computed 
states $u^{(i)}_0(x)$ known over a dense spatial grid, whereas the unlabeled data is the new initial condition $u_0(x_j)$ known only at the sensor locations $x_j\in\mathcal G_s$.

Classification algorithms require a metric that quantifies the proximity of the unlabeled data to the labeled ones. 
The choice of this metric plays a critical role in accurate and consistent classification of data~\cite{Coifman2005,Weinberger,Xing2003}. 
Although the usual $L^2$ norm is often used in this context, there is no theoretical basis for its effectiveness in such classification problems. 
Certain PDEs have an `intrinsic' norm which is dictated by the structure of the governing equations. For instance, the 
intrinsic norm for geodesic equations is induced by the Riemannian metric that defines the geodesic distance. 
In section~\ref{sec:RD Equation}, we show that reaction-diffusion systems also have an intrinsic norm induced by the gradient flow formulation of the equation. 
However, in general, the PDE may not have an immediately obvious intrinsic norm. 
For instance, to the best of our knowledge, the Navier--Stokes equations for fluid flows do not have a preferred metric structure.

Furthermore, the sparsity of the unlabeled data $u_0(x_j)$ introduces an additional difficulty in quantifying its proximity to labeled data.
Most norms used in functional analysis require the functions to be known throughout the spatial domain $\Omega$, or at least on a dense grid over it.
As such, they cannot be computed from sparse measurements.
To address these issues, we introduce a sparsity-promoting metric-learning approach in section~\ref{sec:opt} that
learns a quantifiable metric directly from the labeled data.

\section{Sparsity-promoting metric learning}\label{sec:opt}
In this section, we describe our method for finding the optimal metric which is induced by a norm.
There are many different ways to define a norm on a function space. Here we restrict our search to a class of weighted $L^2$ norms.
For a positive-definite density $\phi:\Omega \to \R^+$, we consider the weighted inner product,
\begin{equation}
    \left<u,v \right>_\phi = \int_\Omega u(x)\cdot v(x) \phi(x)\,\id x, 
\end{equation} 
and the associated norm,
\begin{equation}\label{eq:phi_norm}
    \|u\|^2_\phi = \left<u,u \right>_\phi = \int_\Omega |u(x)|^2 \phi(x) \, \id x.
\end{equation}
This norm induces a metric, where the distance between two states $u$ and $v$ is measured by $\|u-v\|_\phi$. 

We seek a density $\phi$ such that the induced metric has the following properties: (i) It is consistent with the labeled data (as defined shortly below), 
and (ii) It is computable from sparse measurements. To elaborate on the consistency of the metric with the labeled data, we define the \emph{similar set},
\begin{equation}
\mathcal S = \left\{ (u^{(i)}_0,u^{(j)}_0)\ \big|\ \ell(i)=\ell(j)\right\},
\end{equation}
containing the pairs of initial conditions which have the same label. In other words, if $(u^{(i)}_0,u^{(j)}_0)\in \mathcal S$, the initial conditions
$u^{(i)}_0$ and $u^{(j)}_0$ converge to the same attractor as time $t$ tends to infinity.
We also define the \emph{dissimilar set}, 
\begin{equation}
\mathcal D = \left\{ (u^{(i)}_0,u^{(j)}_0)\ \big|\ \ell(i)\neq\ell(j)\right\},
\end{equation}
which contains the pairs of initial conditions that asymptotically converge to different attractors.
We denote the cardinality of similar and dissimilar sets by $|\mathcal S|$ and 
$|\mathcal D|$, respectively. If there are $n_\ell$ labeled data per attractor, then $|\mathcal S| = k{n_\ell\choose 2}$ and $|\mathcal D| = n_\ell^2k (k-1)/2$,
where $k$ is the number of distinct attractors $\mathcal A_\beta$.

We say a metric is consistent with the labeled data if the pairs of initial conditions in the similar set are relatively close as measured in this metric, 
while the pairs in the dissimilar set are relatively far from each other. To quantify these statements, we define the \emph{similar sum},
\begin{equation}
S(\phi) = \sum_{(u_0^{(i)},u_0^{(j)})\in \mathcal S} \| u_0^{(i)}  - u_0^{(j)}\|_\phi^2,
\end{equation}
which measures the sum of pairwise squared distances between the initial states in the similar set $\mathcal S$.
Similarly, we define the \emph{dissimilar sum},
\begin{equation}\label{eq:dissum}
D(\phi) = \sum_{(u_0^{(i)},u_0^{(j)})\in \mathcal D} \| u_0^{(i)}  - u_0^{(j)}\|_\phi^2,
\end{equation}
which measures the sum of pairwise squared distances between the initial states in the dissimilar set $\mathcal S$.
{\cb  Note that the similar sum and the dissimilar sum depend linearly on the density $\phi.$}
The metric $\|\cdot\|_\phi$ is consistent with the data if the similar sum $S(\phi)$ is significantly smaller than 
the dissimilar sum $D(\phi)$.

Xing et al.~\cite{Xing2003} formulated an optimization problem to learn such a consistent metric for finite-dimensional data.
We take a similar approach here, but we generalize their formulation to the infinite-dimensional function space. Furthermore,
we add a sparsity-promoting elastic-net penalization to the optimization problem to ensure that the metric is computable from sparse measurements.
As we show in Appendix~\ref{app:rank1}, the formulation of Xing et al.~\cite{Xing2003} generally returns a degenerate density $\phi$. 
The penalization term that we add to the optimization problem serves the secondary purpose of ruling out these degenerate solutions.

The optimization problem seeks to learn the density $\phi $ by minimizing the sum of pairwise distances between 
initial conditions in the similar set $\mathcal S$. A constraint is added to this optimization to ensure that the distance between
dissimilar initial conditions is above a threshold. More precisely, we solve the constrained optimization problem,
\begin{subequations}\label{original_opt}
\begin{equation} \label{first}
    \inf_{\phi\in L^2(\Omega)}\Big[S(\phi)+ \a \big(\lambda\|\phi\|_1+(1-\lambda)\|\phi\|_2  \big)  \Big],
\end{equation}
\begin{equation} \label{second}
D(\phi)\geq 1,
\end{equation}
\begin{equation}\label{third}
\phi > 0, 
\end{equation}
\end{subequations}
where $\a \left(\lambda\|\phi\|_1+(1-\lambda)\|\phi\|_2  \right)$ is the elastic-net regularization, with sparsity parameter $\lambda\in[0,1)$ and regularization weight $\a \geq 0$.

Lets first consider the case with no regularization, i.e., $\a=0$. 
In this case, equation~\eqref{first} minimizes the similar sum $S(\phi)$, i.e., the total squared distance between labeled initial conditions that converge to the same attractor.
With no further constraint, the density $\phi$ can be chosen arbitrarily small, leading to the trivial solution $\phi=0$ in the limit. 
Constraint~\eqref{second} rules out this trivial solution.
It also ensures that differently labeled data are relatively far from each other by requiring the dissimilar sum $D(\phi)$ to be larger than the threshold $1$.
As discussed shortly, the value of this positive threshold is irrelevant because of the scale invariance of the optimization problem.
Finally, constraint~\eqref{third} ensures that the density $\phi$ is positive definite, and therefore expression~\eqref{eq:phi_norm} is in fact a norm.

The following theorem states some important properties of the optimization problem~\eqref{original_opt}, including its scale invariance and the importance of the elastic-net regularization.
\begin{theorem}\label{thm:opt}
The optimization problem~\eqref{original_opt} has the following properties.
\begin{enumerate}
\item It is convex.
\item The minimizer satisfies $D(\phi)=1$.
\item Scale invariance: If the constraint~\eqref{second} is replaced by $D(\phi)\geq c$, for some $c>0$, then the minimizer changes to $c\phi$.
\item If $\alpha=0$ or $\lambda=1$, then the minimizer is generically non-zero only at a single point $x\in\Omega$.
\end{enumerate}	
\end{theorem}
\begin{proof}
See Appendix~\ref{app:rank1}.
\end{proof}

The convexity of the problem ensures that the optimization problem has a unique minimizer.
The second property in Theorem~\ref{thm:opt} states that the minimum is attained on the boundary of the constraint~\eqref{second}.
This property can be used as a validation test for numerically obtained solutions of the optimization problem.
The third property, i.e., scale invariance, implies that the threshold in constraint~\eqref{second} can be an arbitrary positive number.
Requiring $D(\phi)\geq c$, for some positive constant $c$, would simply change the optimizer $\phi$ to $c\phi$.
As a result, all pairwise distances $\|\cdot\|_\phi$ between data points would be uniformly multiplied by a factor $c$, leaving the subsequent steps (i.e., classification) unaffected.
The last property in Theorem~\ref{thm:opt} highlights the importance of the $L^2$ regularization term $\|\phi\|_2$.
Without this regularization, the optimization problem returns a degenerate rank-1 solution which is too restrictive for our classification purposes.

We also introduce a sparsity-promoting penalization  $\|\phi\|_1$ to the optimization problem to ensure that the resulting norm $\|\cdot\|_\phi$ is computable from sparse measurements.
If $\lambda=0$, the optimal density $\phi$ may be nonzero throughout the domain $\Omega$, and consequently the induced norm $\|\cdot\|_\phi$ would not
be computable from sparse measurements.
As the sparsity parameter $\lambda$ increases, the density $\phi$ tends towards zero throughout the domain, except at points where the state values are critical for minimizing the similar sum $S(\phi)$ while satisfying the constraints. 
This introduces a trade-off between minimizing the total pairwise distance $S(\phi)$ between similarly labeled data and the computability of the norm $\|\cdot\|_\phi$ from sparse measurements. 
For a sparsity parameter $\lambda$, we denote the solution to the optimization problem~\eqref{original_opt} 
by $\phi_\lambda$.

We note that, strictly speaking, to promote sparsity one would use an $L^0$ penalization $\|\phi\|_0 := \mbox{Vol}(\mbox{supp}\,\phi)$, where $\mbox{supp}\,\phi$ denotes the support of $\phi$ and $\mbox{Vol}$ denotes the volume of a set. 
The $L^0$ penalization minimizes the subset of the spatial domain $\Omega$ where $\phi$ is nonzero.
It is well-known that, even in finite dimensions, the optimization problem with the $L^0$ penalization is nonconvex and computationally intractable~\cite{Tibshirani1996}. 
Therefore, we adopt the common compromise by replacing the $L^0$ penalization with $L^1$ penalization which makes the problem convex and computationally tractable, while still 
yielding sparse solutions~\cite{Callaham2019,Ng2004,Tibshirani1996,Hastie2005}.

Finally, we comment on the choice of the regularization weight $\a>0$. This parameter should be large enough so that the sparsity of the solution is ensured
through the penalization term $\a\lambda\|\phi\|_1$. At the same time, it should be small enough to ensure that the regularization term does not overwhelm the similar sum $S(\phi)$. As such, the regularization parameter $\a$ and the sparsity parameter $\lambda$ are free parameters that should be adjusted based on the application.
{  \cb To this end, homotopy continuation algorithms can be used to efficiently determine the suitable free parameters.~\cite{Malioutov2005,Donoho2008,Bieker2020}.}

{\cb Such efficient algorithms would specially be needed for two- and three-dimensional domains. To solve 
the optimization problem~\eqref{original_opt}, the labeled data is discretized over the spatial domain $\Omega$. 
The computational cost of the optimization problem grows with the size of the spatial grid.
} 

We emphasize that the optimization problem~\eqref{original_opt} only uses the labeled data $u_0^{(i)}$ from the precomputed library.
Once an optimal sparse density $\phi_\lambda$ is obtained, it is used to classify new initial conditions $u_0$ which are only known at the sparse locations $x_j\in \mathcal G_s$.
As we show in section~\ref{sec:Numerical Results}, the sparse optimal norm also informs the sensor placement, i.e., the points $x_j$ where state measurements should be gathered.

{\cb
Although here we only consider weighted $L^2$ norms~\eqref{eq:phi_norm}, other weighted $L^p$ or Sobolev norms~\cite{Ukhlov2009} can be used
in a similar fashion. In the case of Sobolev norms, however, clusters of closely placed sensors may be needed to approximate the spatial derivatives.
}

\section{Numerical results}\label{sec:Numerical Results}
In this section, we apply the sparsity-promoting metric-learning optimization to a reaction-diffusion equation {\cb and the FitzHugh--Nagumo equation}. 
In section~\ref{sec:RD Equation}, we first review the reaction-diffusion equation in a general setting and show that there is an intrinsic norm associated with this PDE. 
In section~\ref{sec:RDE_prediction}, we present the numerical results for a one-dimensional reaction-diffusion equation with four asymptotically stable steady states.
{\cb In section~\ref{sec:FHN}, we present the results for the FitzHugh--Nagumo equation with two asymptotically stable steady states. }

\subsection{Reaction Diffusion Equation}\label{sec:RDE}

\subsubsection{Reaction-diffusion equation: intrinsic norm} \label{sec:RD Equation}
We consider the reaction-diffusion equation,
\begin{equation}\label{generic_RDE}
\frac{\d u}{\d t}  =  \nu\frac{1}{w} \grad \cdot \big( w \grad u  \big)  + f(u),\quad  u(x,t_0)=u_0(x),
\end{equation}
with the Neumann boundary conditions $n\cdot \grad u  = 0$. Here, $u:\Omega\times \mathbb R^+\to \mathbb R$ denotes the system state which belongs to a Hilbert function 
space $H$ to be specified, and $w:\Omega\to \mathbb R^+$ is a prescribed weight. The reaction term $f(u)$ is a nonlinear function of the state $u$, and $\nu>0$ is the diffusion coefficient.
Note that $w\equiv 1$ corresponds to the usual homogeneous diffusion $\nu \Delta u$. Here, we consider the more general case where the weight $w$ is not constant everywhere.

It is not immediately obvious what norm is suitable 
for measuring the distance between two states of the reaction-diffusion equation~\eqref{generic_RDE}. 
While the $L^2$ distance is commonly used, this PDE has an intrinsic norm
dictated by the gradient flow structure of the equations. To see this structure, consider 
the energy functional,
\begin{equation} \label{energy_functional}
    I[u] = \int_\Omega \left[ \frac{\nu}{2}|\grad u|^2 - F(u) \right] w \, \id x,
\end{equation} 
where $F(u)$ satisfies $F'(u) = f(u)$. 
Next we consider the gradient flow,
\begin{equation}\label{eq:gradflow}
\frac{\d u }{\d t}  = -\grad_H I(u),
\end{equation}
associated with the energy functional~\eqref{energy_functional}.
The gradient of $I$ at $u\in H$ is an element $\nabla_H I(u)\in H$ that 99satisfies,
\begin{equation}
    \left<\grad_H I(u) , v \right>_H = \frac{\id }{\id \varepsilon } I[u+\varepsilon v]\Big|_{\varepsilon=0} \text{, } \quad \forall v\in H. 
\end{equation}
Here, $\left<\cdot,\cdot \right>_H$ denotes the inner product in $H$.
If we use the $L^2$ inner product, i.e. $H = L^2(\Omega)$, the gradient flow~\eqref{eq:gradflow} would be a PDE that differs 
from the reaction-diffusion equation~\eqref{generic_RDE}. As we show in Appendix~\ref{app:inorm},
to retrieve this reaction-diffusion equation, we must use the weighted inner product,
\begin{equation} \label{eqn:weight_product}
\left<u,v \right>_w := \int_\Omega u(x)\cdot v(x)w(x)\, \id x.
\end{equation}
In this case, $-\nabla_H I(u)$ coincides with the right-hand side of the reaction-diffusion equation~\eqref{generic_RDE}.
We refer to the weighted inner product $\langle \cdot,\cdot\rangle_H = \langle \cdot,\cdot\rangle_w$ as the intrinsic inner product of
equation~\eqref{generic_RDE} since it is the only inner product consistent with its gradient flow structure.
The inner product $\langle \cdot,\cdot\rangle_w$  induces the intrinsic norm,
\begin{equation}\label{eq:wnorm}
\| u \|_w^2 := \int_\Omega |u(x)|^2w(x)\id x.
\end{equation}
In the special case, where $w\equiv 1$, the intrinsic norm coincides with the $L^2$ inner product.

{\cb We point out that other PDEs, which can be interpreted as gradient flows, 
often have an intrinsic metric associated with them. Examples include the
porous medium equation~\cite{Otto2001} and the Fokker--Planck equation~\cite{Otto1998}.
We refer the reader to Ref.~\cite{ambrosio2008} for an exhaustive review.}

In section~\ref{sec:RDE_prediction}, we consider a specific reaction-diffusion equation of the form~\eqref{generic_RDE} which has four asymptotically stable steady states.
We use a precomputed library of initial conditions, labeled according to their asymptotic behavior, to predict the asymptotic state of unlabeled initial conditions.
We compare the classification results when different norms are used to quantify the proximity of the unlabeled data to the precomputed library. 

\subsubsection{Prediction Results : Reaction Diffusion Equation} \label{sec:RDE_prediction}
In this section, we apply the SPML optimization~\eqref{original_opt} to a one-dimensional reaction-diffusion equation~\eqref{generic_RDE} with the reaction term
$f(u) = -u \left(\frac{1}{2}-u\right)(1-u)$. We consider the spatial domain $\Omega = [-1,1]$, so that $ u:[-1, 1]\times \R^+\to \R$.
The resulting PDE reads,
\begin{equation}
    \frac{\d u}{\d t}=  \nu\frac{1}{w}\frac{\d}{\d x}\left(w \frac{\d u }{\d x}\right) - u \left(\frac{1}{2}-u\right)(1-u),
    \label{eq:RDE-1D}
\end{equation}
with the Neumann boundary conditions $\partial_x u|_{x=\pm 1}=0$, and the diffusion coefficient $\nu =10^{-2}$.
We discretize this equation using $201$ equispaced points in space $\Omega=[-1,1]$, resulting in
the grid size $\Delta x=0.01$. We discretize the integrals in the optimization problem~\eqref{original_opt} using the same grid and
the trapezoidal rule. The reaction-diffusion equation is numerically integrated in time by an embedded Runge--Kutta scheme~\cite{ode45} (MATLAB's \texttt{ode{45}}), 
with relative and absolute tolerance of $10^{-5}$. 
\begin{figure}
	\centering
	\subfigure[\label{fig:weight}]
	{\includegraphics[width=0.48\textwidth]{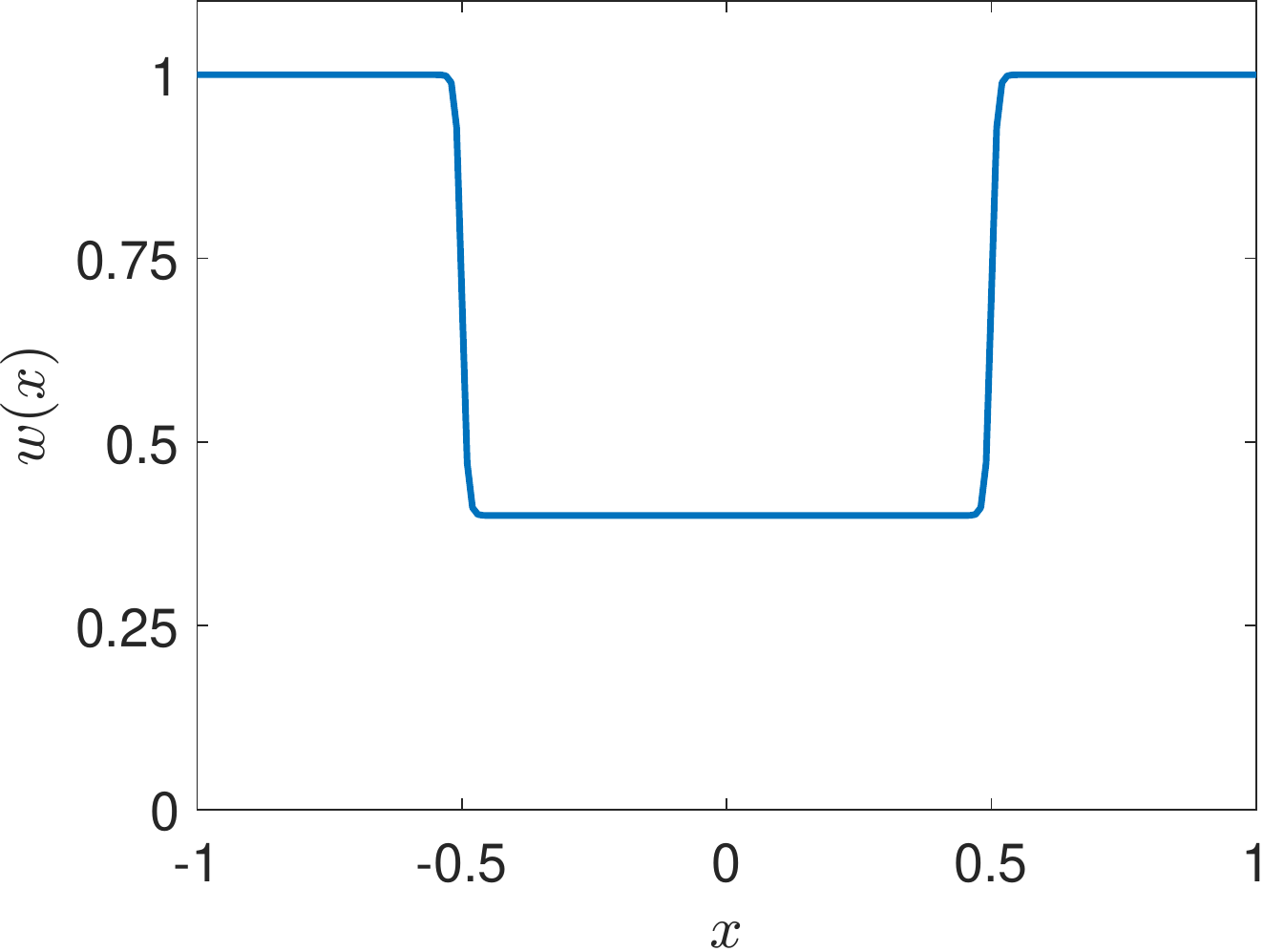}}
	\subfigure[ \label{fig:steady_states}]
	{ \includegraphics[width=0.48\textwidth]{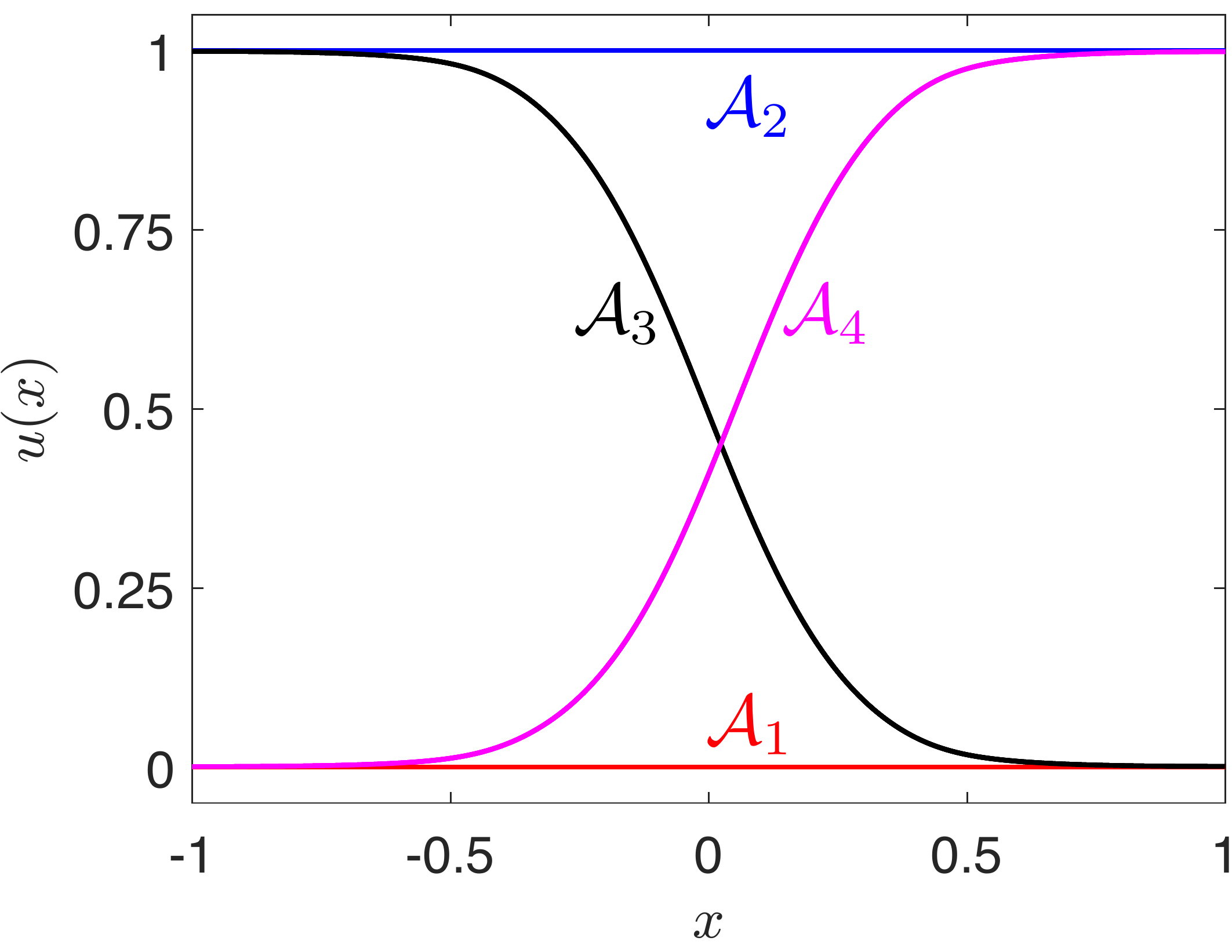}}
	\caption{(a) Weight function~\eqref{eq:w}. 
		(b) Stable steady states of reaction-diffusion equation~\eqref{eq:RDE-1D}.} 
\end{figure}

With the constant weight $w\equiv 1$, the reaction-diffusion equation~\eqref{eq:RDE-1D} has three steady states at $u=0$, $1/2$ and $1$.
The steady state $u=1/2$ is unstable, and therefore the system has only two asymptotically stable attractors at $u=0$ and $u=1$.
To make the classification problem more challenging, we instead use the weight function, 
\begin{equation}\label{eq:w}
    w(x) = a\tanh\left(\frac{x-x_0}{\epsilon}\right)+a\tanh\left(\frac{-x-x_0}{\epsilon}\right)+1,
\end{equation}
where $a=0.3$, $x_0=0.5$, and $\epsilon=0.01$. This function, shown in figure~\ref{fig:weight},
is approximately $0.4$ in the middle section of the domain and converges towards $1$ near the boundaries, with sharp transitions at $x=0.5,\, x=-0.5$. 
The weight function reduces the amount of diffusion that occurs in the middle of the domain allowing the reaction term to dominate diffusion in that part of the domain.
As a result, the reaction-diffusion equation with the weight~\eqref{eq:w} has four asymptotically stable attractors as shown in figure~\ref{fig:steady_states}.
The constant steady states $u\equiv 0$ and $u\equiv 1$ can be shown to be stable using linear stability analysis. 
The other two steady states are found empirically from the numerical evolution of thousands of initial conditions. 

Given sparse measurement of an initial condition $u_0$, we would like to predict to which of the four attractors it converges. 
To this end, we use a library of precomputed initial conditions which are given on a dense grid over $\Omega$. 
These initial conditions are evolved under the reaction-diffusion equation and labeled according to the attractor that they converge to. 
Given sparse measurements of an out-of-library state, we use a clustering algorithm to classify this state, i.e., predict its asymptotic attractor. 
We summarize these steps in Algorithm~\ref{alg:SPML} and expand on each step below.
   
\begin{algorithm}[H]
     \caption{Classification based on sparsity-promoting metric-learning (SPML)}
     \label{alg:SPML}
     \begin{algorithmic}[1]
			\State Generate initial conditions $u_0^{(i)}(x)$ on a dense grid.
			\State Evolve each initial condition in time and assign a label $\ell(i) \in\{1,2,3,4\}$ to it based on the attractor it converges to.
			\State Solve the optimization problem~\eqref{original_opt} to obtain the optimal density $\phi_\lambda$.
			\State Measure a new initial condition $u_0(x_j)$ on a sparse grid, $x_j\in \mathcal G_s$.
            \State Classify the new initial conditions using the precomputed library $u_0^{(i)}$ {\cb using $\|\cdot\|_{\phi_\lambda}$:}  
            \begin{algsubstates}
            	\State {\cb Find the closest labeled state to the new initial condition $u_0$: $\bar i=\arg\min_i \|u_0^{(i)}-u_0\|_{\phi_\lambda}$}
            	\State {\cb Predict that $u_0$ converges to the steady state $\ell(\bar i)$.}
            \end{algsubstates}
\end{algorithmic}
\end{algorithm}

\begin{figure}
	\centering
	\subfigure[\label{fig:initialconditions_1}]{\includegraphics[width=0.49\textwidth]{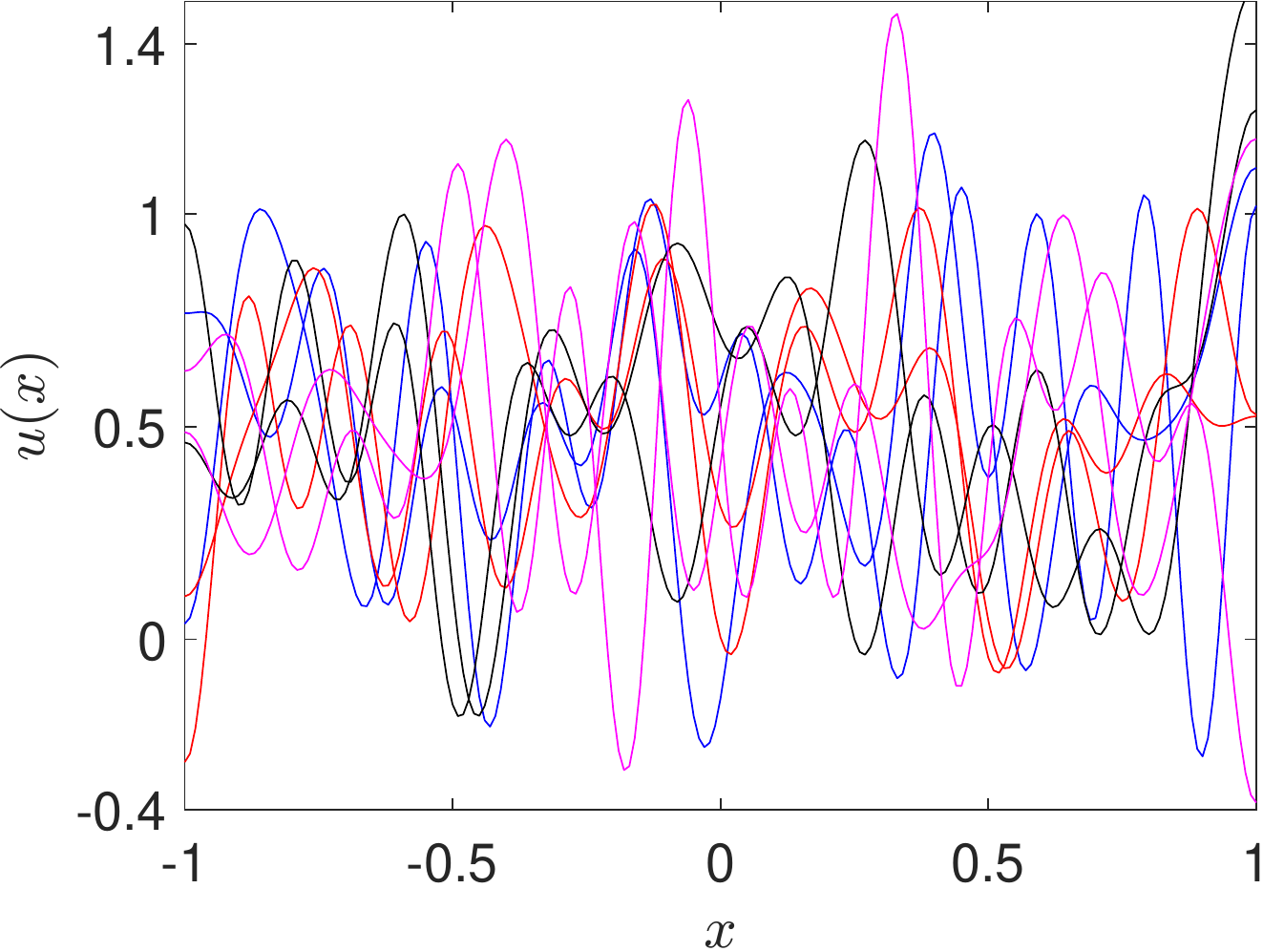}}
    \subfigure[\label{fig:initialconditions_2}]{\includegraphics[width=0.49\textwidth]{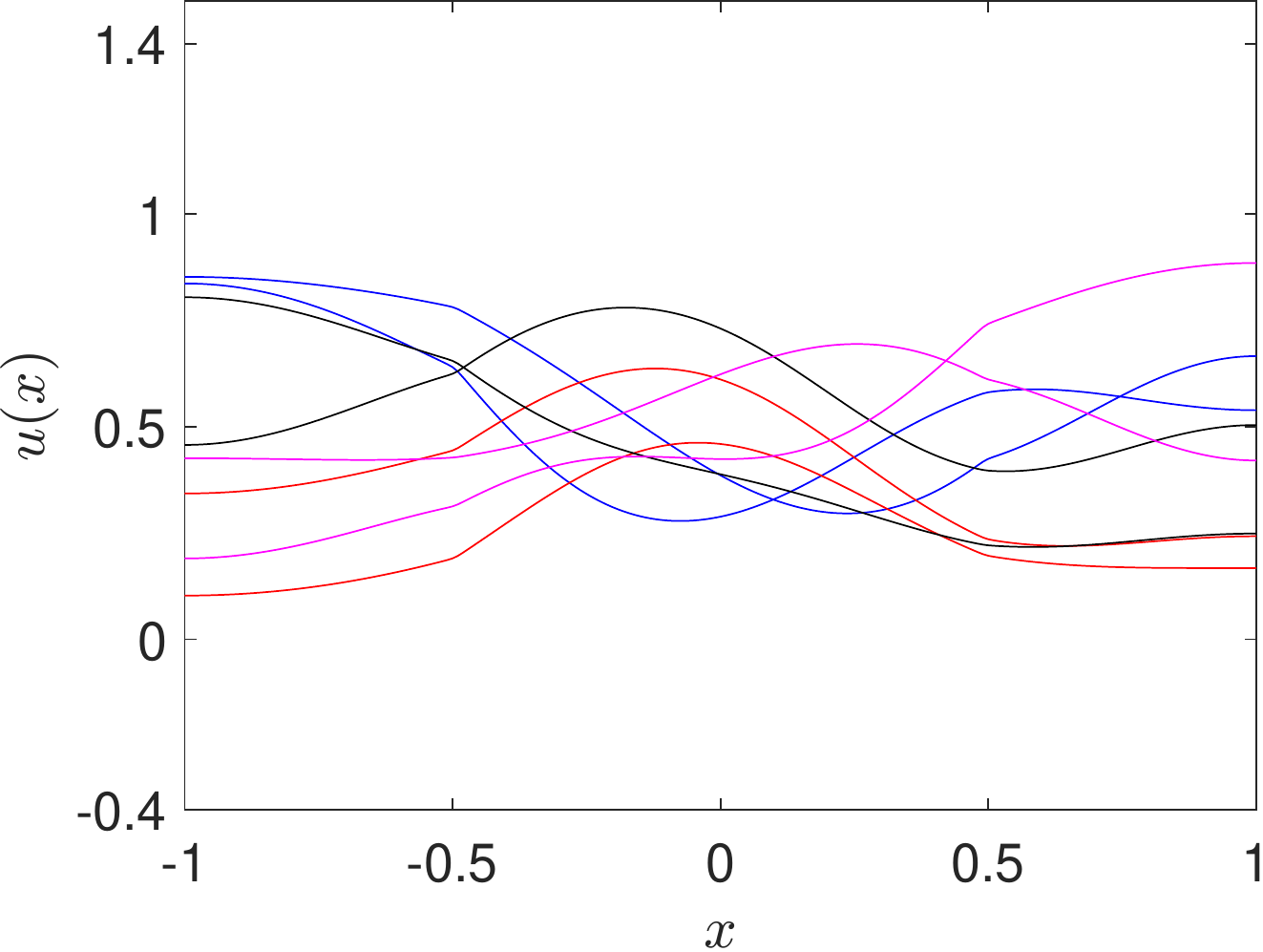}}
    \caption{(a) Eight functions $u(x,0)$ randomly generated based on equation~\eqref{eq:IC}. (b) The time evolution of these functions to time $t_0=10$. Each state is colored according to the attractor to which it asymptotically converges; see figure~\ref{fig:steady_states}.}
	\label{fig:initialconditions}
\end{figure}

In steps 1 and 2 of Algorithm~\ref{alg:SPML}, we generate the initial conditions $u_0^{(i)}$ that constitute the labeled precomputed library.
To this end, we generate random smooth functions,
\begin{align} \label{eq:IC}
    u(x,0) = & \frac12 + \frac{1}{10} \sum_{k=1}^{10} \Big[a_k\cos(k\pi x) \nonumber \\
    &+ b_k \sin\left(\frac{2k-1}{2}\pi x\right)\Big],
\end{align} 
where $a_k$ and $b_k$ are random numbers with a standard normal distribution. Figure~\ref{fig:initialconditions_1} shows eight random functions generated from equation~\eqref{eq:IC}.
Classification based on these random data, which have not undergone any reaction-diffusion dynamics, proves too demanding. 
Instead, we evolve the random functions under the reaction-diffusion equation~\eqref{eq:RDE-1D} for a short period of ten time units to obtain the initial conditions $u(x,t_0)=u_0(x)$, where $t_0=10$. 
The resulting initial conditions are shown in figure~\ref{fig:initialconditions_2}.
We use the evolved states $u_0(x)$ as the labeled data. This is in line with applications, where state observations are made after the state has evolved for a certain period of time
under the system dynamics.

In step 2, we label the initial conditions. To this end, we evolve the initial conditions until they converge to their asymptotic state. If initial condition $u_0^{(i)}$ 
converges to attractor $\mathcal A_\beta$, we assign the label $\ell(i)=\beta$ to that initial condition. 
On average it takes about $100$ time units for the initial conditions to converge to an attractor.
\begin{figure}
	\centering
	\subfigure[\label{fig:5labels} ]{\includegraphics[width=.49\textwidth]{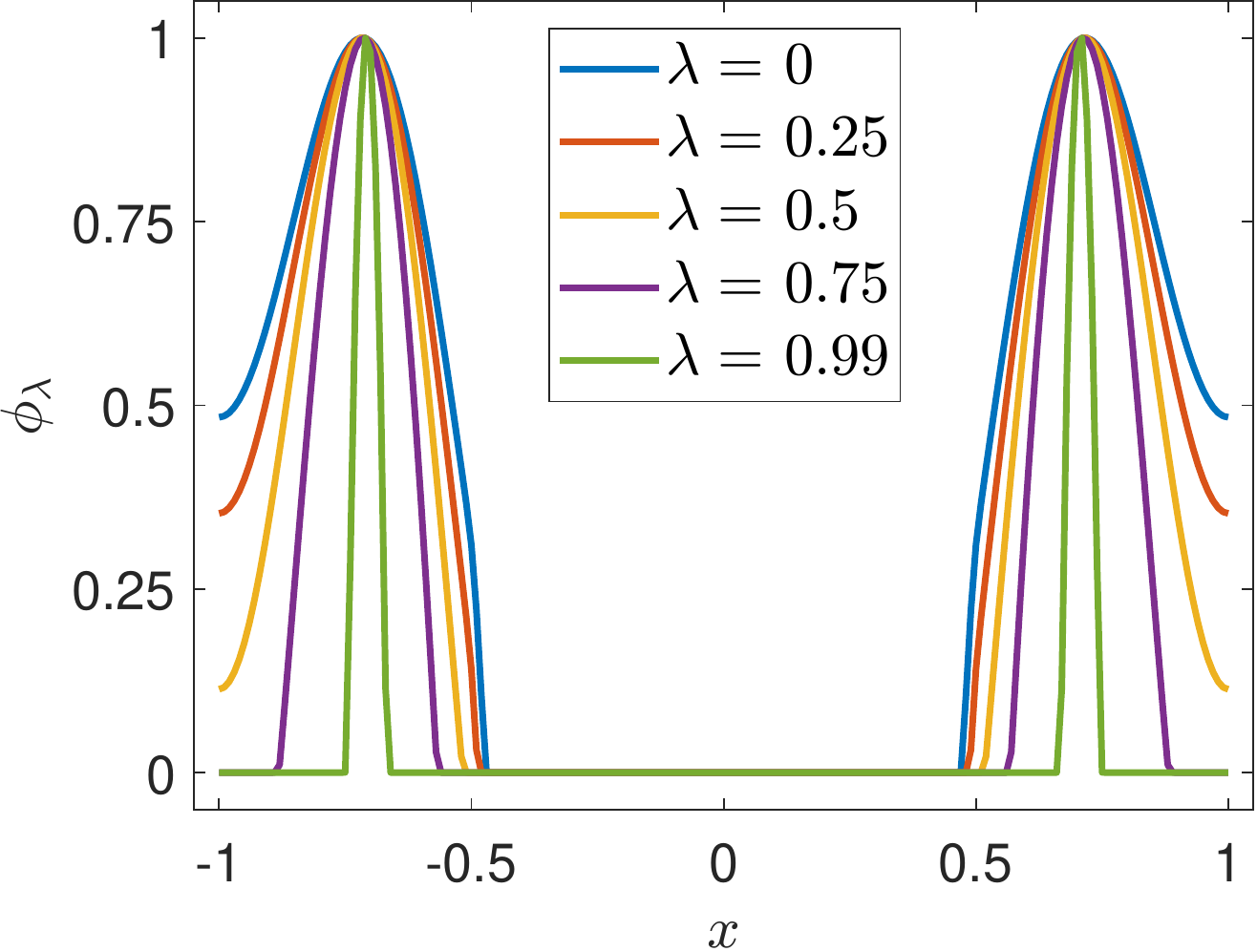}}
	\subfigure[\label{fig:10labels}]{\includegraphics[width=.49\textwidth]{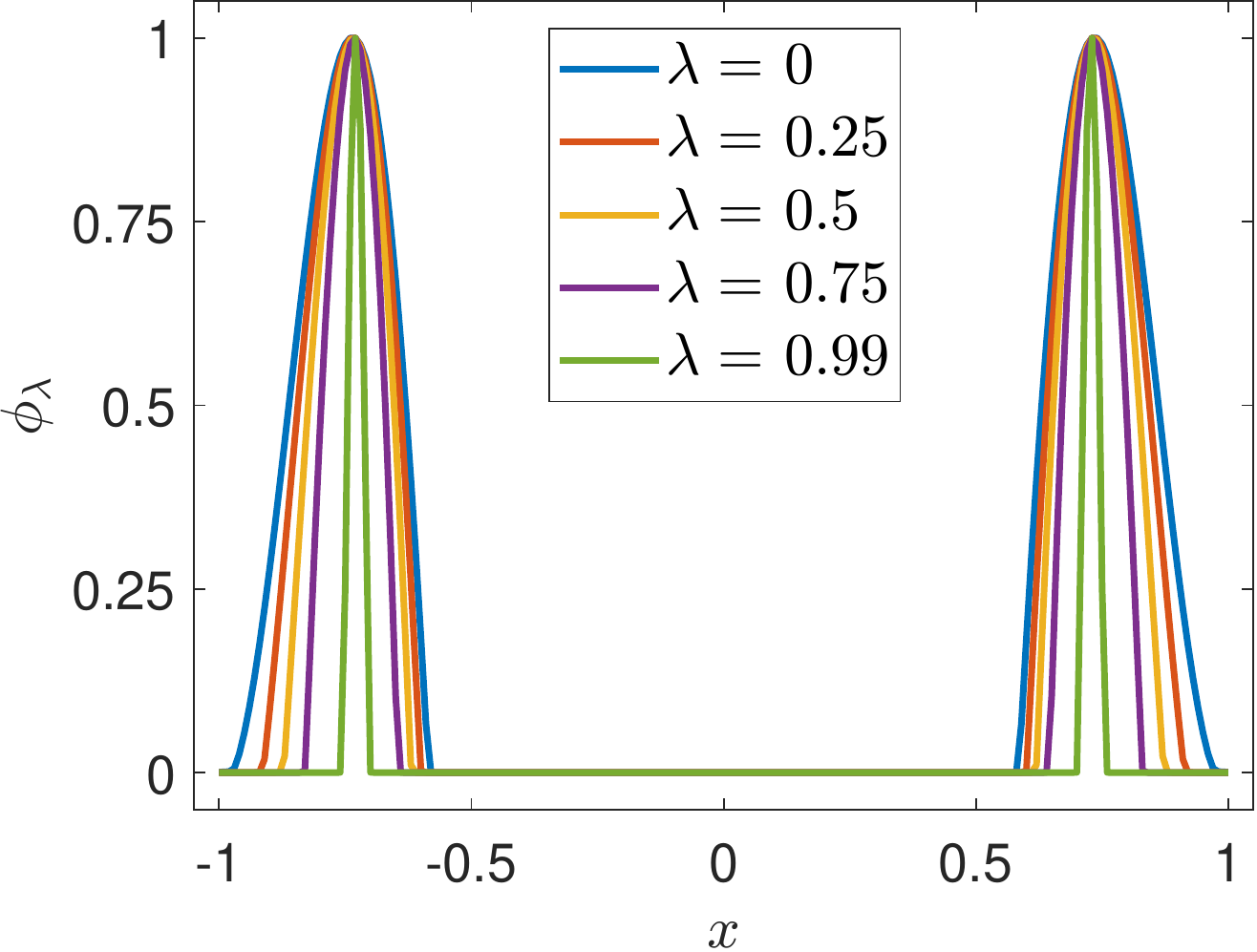}}
	\caption{The learned density $\phi_{\lambda}$ with $\a= 1$ and various values of $\lambda$. 
    The densities are learned from (a) 10 labeled data per attractor. 
    (b) 20 labeled data per attractor. The densities are normalized by their maximum value, $\max_{x\in\Omega} \phi_{\lambda}(x)$.} 
	\label{fig:phi}
\end{figure}

We note that the reaction-diffusion equation~\eqref{eq:RDE-1D} is equivariant under the transformation $x\mapsto -x$.
This is due to the reflection symmetry of the PDE; the change of variable $x\mapsto -x$ leaves the equation unchanged. 
As a result if $u(x,t)$ is a solution to the reaction-diffusion equation corresponding to the initial condition $u_0(x)$, then 
$u(-x,t)$ is the solution corresponding to the initial condition $u_0(-x)$. In particular, this implies that if $u_0(x)$ converges 
to attractor $\mathcal A_1$ (respectively, $\mathcal A_2$), then the initial condition $u_0(-x)$ also converges to $\mathcal A_1$ (respectively, $\mathcal A_2$).
Moreover, if $u_0(x)$ converges to attractor $\mathcal A_3$, then $u_0(-x)$ converges to attractor $\mathcal A_4$ and vice versa. This is due to the fact that $\mathcal A_3$
and $\mathcal A_4$ are reflected copies of each other. Using this equivariance, one can double the size of the precomputed library at no additional computational cost. 
The results reported below include these symmetric copies.



Step 3 constitutes the main part of Algorithm~\ref{alg:SPML} and involves solving the optimization problem~\eqref{original_opt}.
We solve this problem using an interior-point algorithm (Matlab's \texttt{fmincon}) with tolerance $10^{-6}$.
The results reported here correspond to the regularization weight $\alpha=1$, but we have varied $\alpha$ in the interval $[0.5,1]$ and
observed no significant change in the results.
We start with $\lambda=0$, corresponding to the optimization problem with no sparsity promotion. 
Figure~\ref{fig:phi} shows the optimal density $\phi_0$ computed from a library of initial conditions $u_0^{(i)}$ consisting of $10$ and $20$ labeled data per attractor.
This optimal density is non-zero over large portions of the domain $[-1,1]$. 
As such, the corresponding norm (cf. equation~\eqref{eq:phi_norm}) cannot be used to quantify the proximity of states when only sparse measurements are available.
We refer to $\phi_0$, the optimal density at $\lambda=0$, as the \emph{non-sparse optimal density}.

Next, we gradually increase $\lambda$ promoting the sparsity of the optimal density. 
Figure~\ref{fig:phi} shows that as $\lambda$ increases, 
the optimal density $\phi_\lambda$ gradually concentrates near $x=\pm 0.72$ and is zero almost everywhere else. 
As such, the corresponding norm can be quantified from measurements around $x = \pm 0.72$. 

Further increasing the sparsity parameter $\lambda$, the optimal density $\phi_\lambda$ concentrates more and more at $x=\pm 0.72$.
Taking this observation to its logical limit, we introduce the \emph{sparse optimal density} $\phi_s  = \delta (x+0.72) +\delta(x-0.72)$, where $\delta(x)$ is the Dirac delta function.
This density induces the pseudo-metric,
\begin{align}\label{eq:snorm}
    \|u - v\|_{\phi_s}^2 =& \left( u(-0.72)-v(-0.72)\right)^2 \nonumber\\
                          & + \left( u(+0.72)-v(+0.72)\right)^2,
\end{align}
which is computable from two-point measurements. 
As a result, in step 4 of Algorithm~\ref{alg:SPML}, we only make state measurements at $x=\pm 0.72$.
It is worth noting, therefore, that the sparsity-promoting optimization~\eqref{original_opt} informs the optimal sensor placement $\mathcal G_s = \{-0.72,+0.72\}$ in addition to returning an optimal norm for quantifying the proximity of system states.

Finally, step 5 consists of predicting the asymptotic behavior of a new out-of-library initial condition $u_0(x)$.
To this end, we use a clustering algorithm that assigns a label to the new state $u_0$ given the precomputed labeled states $u_0^{(i)}$ (i.e., semi-supervised classification).
Clustering algorithms typically require a norm to quantify the proximity of various states. 
Here we use the sparse optimal norm $\|\cdot\|_{\phi_s}$, which was obtained by solving the SPML optimization~\eqref{original_opt}.
For comparison, we also discuss the prediction results using three other norms: the $L^2$ norm, the intrinsic norm of the reaction-diffusion equation~\eqref{eq:wnorm}, and the {\cb learned} non-sparse optimal norm~\eqref{eq:phi_norm}.
These norms are summarized below:
\begin{enumerate}
	\item The usual $L^2$ norm, $\| \cdot \|_{L^2}$,
	\item The intrinsic norm, $\|\cdot\|_w$ (cf. equation~\eqref{eq:wnorm}),
	\item The non-sparse optimal norm, $\|\cdot \|_{\phi_0}$ (cf. figure~\ref{fig:phi} and equation~\eqref{eq:phi_norm}),
	\item The sparse optimal norm, $\|\cdot\|_{\phi_s}$ (cf. equation~\eqref{eq:snorm}).
\end{enumerate}
We emphasize that only the sparse optimal norm $\|\cdot\|_{\phi_s}$ is computable from sparse measurements, whereas computing the other three norms requires the initial condition to be known on a dense spatial grid. {\cb Also recall that the weight $w$ in the intrinsic norm is dictated by the PDE (see Section~\ref{sec:RD Equation}), whereas the weights $\phi_0$ and $\phi_s$ are obtained from solving the SPML optimization~\eqref{original_opt}.}

\begin{figure}
	\centering
	\includegraphics[width=0.85\textwidth]{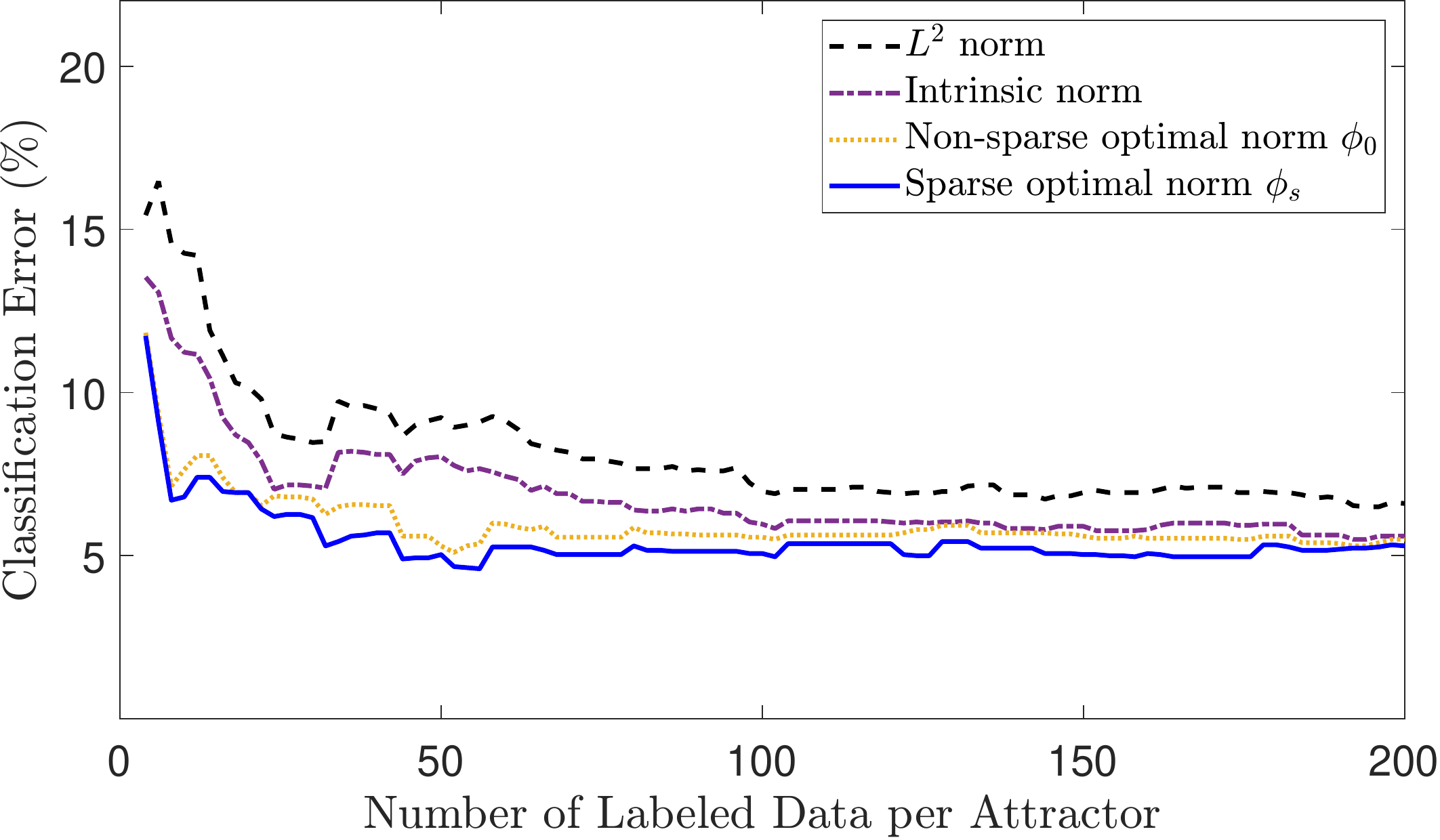}
	\caption{Classification errors as a function of the number of labeled data.}%
	\label{fig:errorsforlabels} 
\end{figure}

Figure~\ref{fig:errorsforlabels} shows the classification errors for various numbers of labeled data.   
In each case, the classification error is computed for 3000 out-of-library initial conditions.
Here, we use nearest neighbor classification, the simplest classification algorithm. 
More precisely, for each out-of-library initial condition $u_0$, we find its closest neighbor from the library $u_0^{(i)}$ of precomputed states. 
The nearest neighbor is determined based on the distance between states as measured by the norms listed above. 
The out-of-library initial condition $u_0$ is then predicted to converge to the same attractor as its nearest neighbor in the precomputed library.

The classification errors corresponding to the non-sparse norms, $\|\cdot\|_{L^2}, \|\cdot\|_{w}$, and $\|\cdot\|_{\phi_0}$ exhibit a similar trend:
as the number of labeled data $n_\ell$ increases, the classification errors generally decrease.
The $L^2$ norm performs poorly as it returns the highest error for any number of labeled data. 
The intrinsic norm $\|\cdot \|_w$ performs slightly better compared to the $L^2$ norm. The learned norm $\|\cdot\|_{\phi_0}$ outperforms both the $L^2$
norm and the intrinsic norm.

Recall that the distances measured in the non-sparse norms, $\|\cdot\|_{L^2}$, $\|\cdot\|_{w},$ and $\|\cdot\|_{\phi_0}$,
are not quantifiable when only sparse measurements of the new initial condition $u_0$ are available. 
To perform classification based on such sparse measurements, we use the sparse optimal norm $\|\cdot\|_{\phi_s}$, which is computed from function evaluations at only two points $x=\pm 0.72$.
The corresponding classification errors, also shown in figure~\ref{fig:errorsforlabels}, are comparable to the classification errors using non-sparse learned norms.
When a small number of labeled data are used, the classification error based on sparse measurements are under $11\%$, allowing correct prediction of asymptotic behavior of more than $89\%$ of the initial conditions. As the number of labeled data increases, the classification error decreases to about $5\%$, resulting in correct predictions $95\%$ of the time.
Often the classification error from sparse learned norm $\|\cdot\|_{\phi_s}$ is lower than the error corresponding to the non-sparse learned norm $\|\cdot\|_{\phi_0}$.
This is perhaps specific to the reaction-diffusion equation~\eqref{eq:RDE-1D}. We do not expect that the classification from sparse measurements would generally outperform classification based on dense measurements.

\begin{figure}
	\centering
	\includegraphics[width=0.8\textwidth]{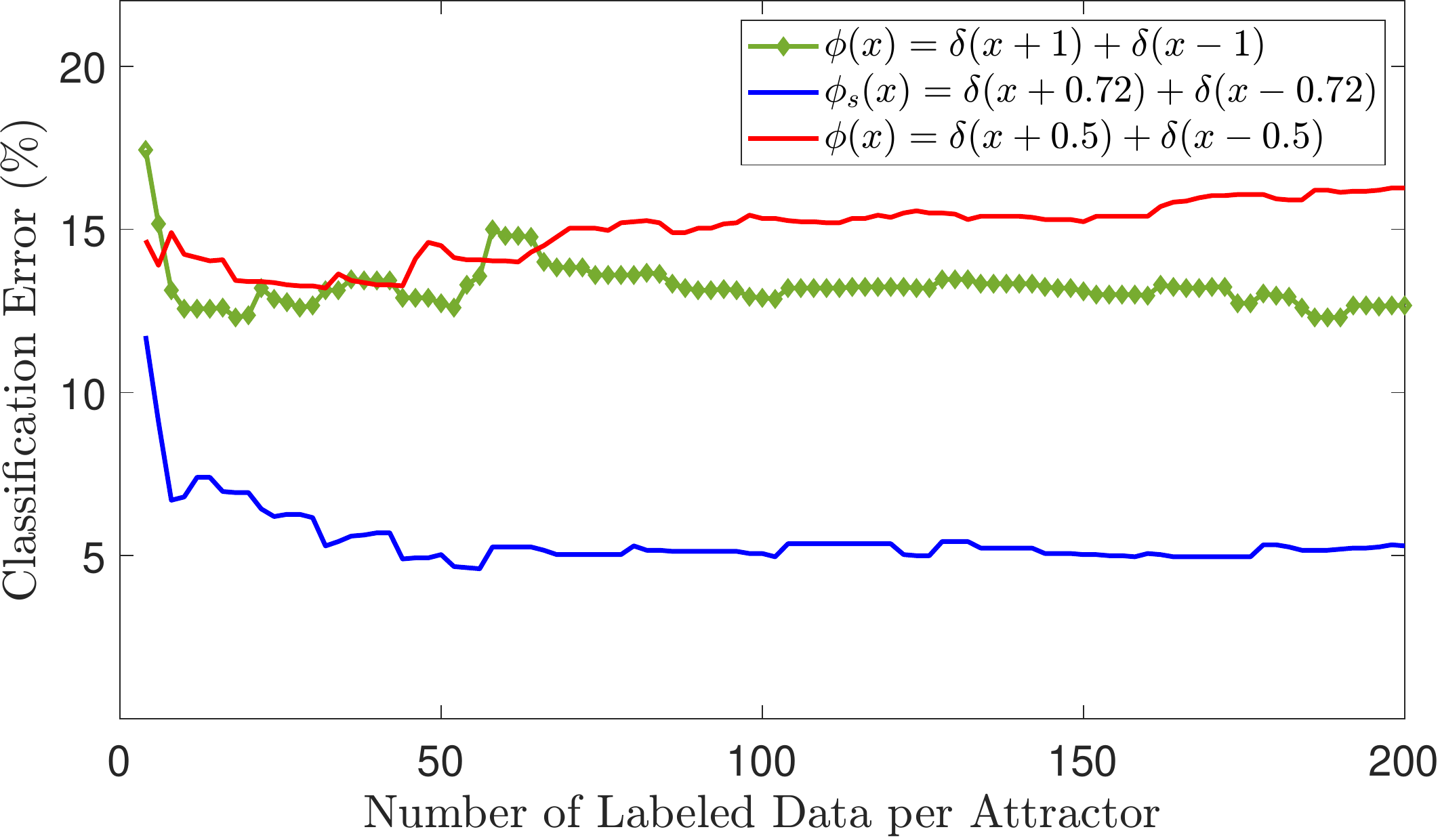}
    \caption{Classification errors for three types of sensor locations. The classifications based on measurements at the optimal points $x=\pm 0.72$ outperform those based on boundary points $x=\pm 1$ or interior points $x=\pm 0.5$.}
	\label{fig:errorsforlabels2}
\end{figure}

Figure~\ref{fig:errorsforlabels2} demonstrates the optimality of the points $x=\pm 0.72$ for sensor placement.
For comparison, it shows the classification error based on two-point measurements if the sensors were to be located at the boundaries $x=\pm 1$ 
or closer to the center of the domain at $x=\pm 0.5$. The distance between states would then be measured based on the densities $\phi(x) = \delta(x+1) + \delta(x-1)$
and $\phi(x) = \delta(x+0.5) + \delta(x-0.5)$. As shown in figure~\ref{fig:errorsforlabels2}, classification based on the learned optimal locations $x=\pm 0.72$
outperforms classifications based on other two-point measurements. 

\begin{figure}
\centering
\includegraphics[width=0.5\textwidth]{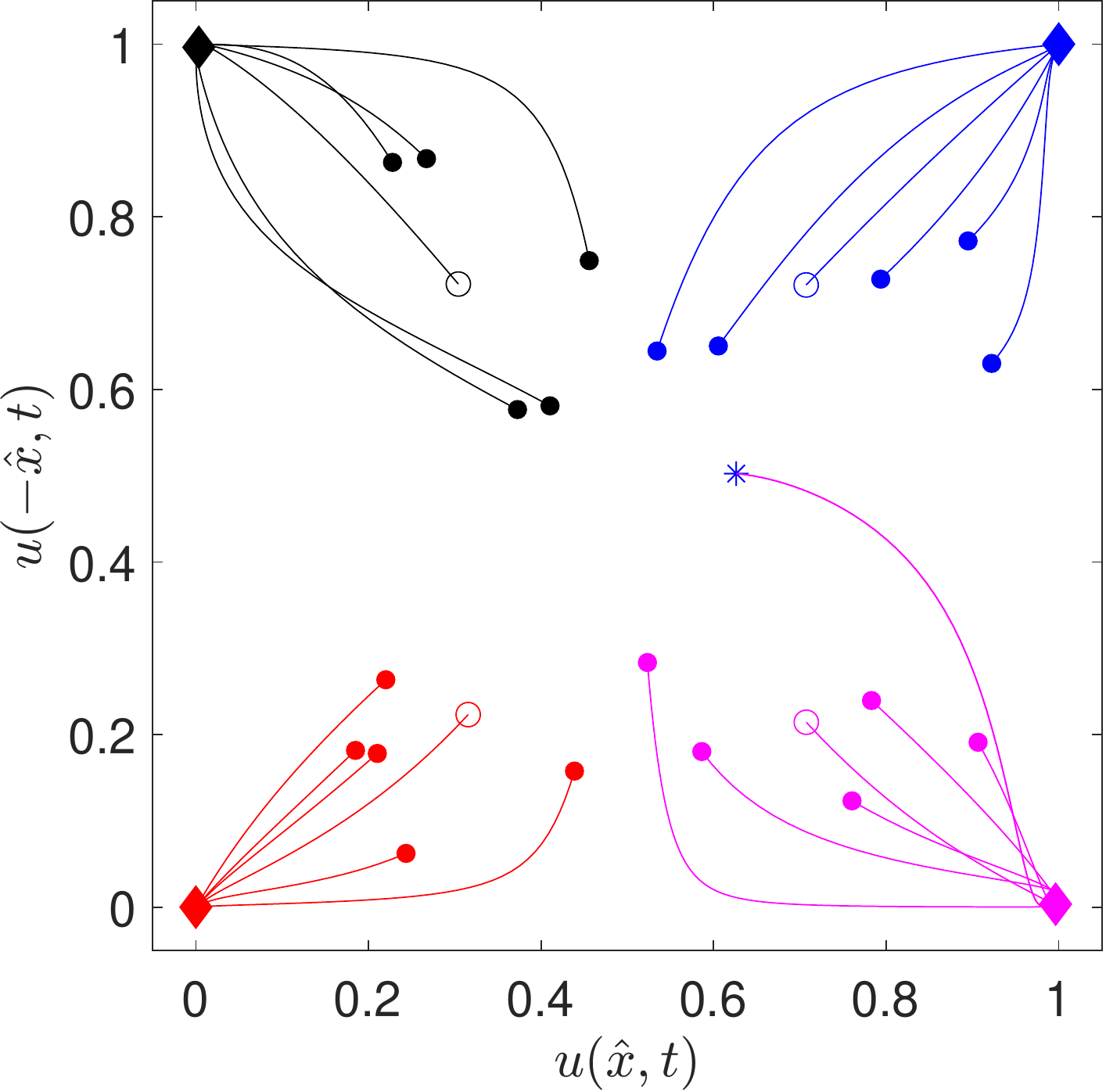} 
\caption{Reaction-diffusion trajectories in the observation space $(u(\hat x,t),u(-\hat x,t))$, where $\hat x=0.72$.
The attractors are marked by the diamond symbols. For each attractor 5 labeled initial conditions $u_0^{(i)}$ (filled circles) 
from the precomputed library and their trajectories (solid curves) are plotted. For each attractor, one unlabeled initial condition (empty circles) and its trajectory are also shown.
{\cb One misclassified initial condition (marked by a star symbol) and its trajectory are also shown.} 
}
\label{fig:projTraj}
\end{figure}

The two-point measurements at the optimal points $x=\pm 0.72$ can be used to produce a low-dimensional visualization
of the system dynamics. To this end, we define the map $\pi:H\to \mathbb R^2$ such that $\pi\circ u = (u(\hat x,t), u(-\hat x,t))$, where $\hat x = 0.72$. As such, the map $\pi$ is a projection from the function space $H$ to the observation space
$\mathbb R^2$. In the observation space, the attractors lie on the corners of the unit square. For instance, the 
steady state $u\equiv 0$ (respectively, $u\equiv 1$) is represented by the point $(0,0)$ (respectively, (1,1)). 
Figure~\ref{fig:projTraj} shows the steady states (diamond symbols) in the observation space.
We also plot {\cb $25$} trajectories 
of the system starting from various initial conditions. 
In the observation space, the initial conditions that converge to different attractors are well-separated, elucidating the success of the nearest neighbor classification. 
{\cb In figure~\ref{fig:projTraj}, we also show one misclassified initial condition (marked by a star symbol). This initial condition
is closer to the blue labeled data and therefore is predicted to converge to attractor $\mathcal A_2$. 
However, in reality, it converges to attractor $\mathcal A_4$. Nonetheless, as figure~\ref{fig:errorsforlabels} shows,
the number of such misclassified data points is very low.}

We point out that other clustering algorithms can be employed to perform the classification. Most clustering algorithms also need a measure of proximity between data points.  For instance, spectral clustering requires forming an affinity matrix $W$ whose elements $w_{ij}$ measure the proximity of pairs of data points~\cite{Bertozzi2018}. Our SPML optimization~\eqref{original_opt}
can be employed to identify an optimal measure of proximity to be used in such clustering algorithms.

\subsection{FitzHugh--Nagumo Equation} \label{sec:FHN}
The FitzHugh--Nagumo equation is a slow-fast dynamical system consisting of two coupled PDEs. 
This equation was first proposed as a model for propagation of electrical signals in neurons~\cite{fitzhugh1955,Nagumo1962}.
FitzHugh--Nagumo model consists of two spatiotemporal variables $u(x,t)$ and $v(x,t)$ which satisfy
\begin{align}\label{FHN_eqn}
    \frac{\d u}{\d t} & = \nu \Delta u - v +f(u)\nonumber\\
    \frac{\d v}{\d t} & = \beta u - \gamma v ,
\end{align}
where $f(u)$ is a nonlinear reaction term, and $u,v: \Omega\times\R^+\to \R$ denote the state variables.
Here we consider the one-dimensional case with the spatial domain $\Omega = [-1,1]$, the reaction term
$f(u) = u \left(\frac{1}{2}-u\right)(u-1)$, and Neumann boundary conditions  $n\cdot \grad u  = 0$ and  $n\cdot \grad v=0$. 
The model parameters are $\beta = 10^{-2}$, $\gamma = 1$, and the diffusion coefficient $\nu=10^{-2}$.

We discretize this equation using $201$ equispaced points in the domain $\Omega=[-1,1]$, resulting in
the grid size $\Delta x=0.01$.
We discretize the integrals in the optimization problem~\eqref{original_opt} using the same grid and the trapezoidal rule. 
The FitzHugh--Nagumo equation is numerically integrated in time by an embedded Runge--Kutta scheme~\cite{ode45}, 
with relative and absolute tolerance of $10^{-5}$. 

Spatially constant steady state solutions of the FitzHugh--Nagumo equation can be obtained by a straightforward calculation.
Assuming the steady states to be constant in space, $u(x,t)\equiv u$ and $v(x,t)\equiv v$, they satisfy
\begin{align*}
    0 & = -v +u\left(\frac{1}{2}-u\right)(u-1) \\
    0 & ={\beta}u -\gamma v
.\end{align*}
Solving for $u$, we obtain
 $\frac{\b}{\g} u =u\left(\frac{1}{2}-u\right)(u-1)$ which admits the solutions,
\begin{equation}
u = 0,\quad \frac{3 \g \pm \sqrt{\g^2-16\b\g}}{4\g}.
\end{equation}

Using $\beta = 10^{-2}$ and $\gamma = 1$, the steady state solutions for $u$ are $u_1=0$, $u_2 \approx 0.5209$, or $ u_3\approx 0.9791$,
with the corresponding $v_ i = \beta u_i/\gamma$.
Linear perturbation analysis of the FitzHugh--Nagumo equation shows that $u_1$ and  $u_3$ are asymptotically stable attractors while $u_2$ is an unstable steady state.

We apply the SPML algorithm~\ref{alg:SPML} to the FitzHugh--Nagumo equation to predict whether an initial condition converges to the stable steady state $(u_1,v_1)$ or $(u_3,v_3)$.
To make the prediction task somewhat more challenging, we assume that only variable $u$ is observable.
In other words, although variable $v$ contributes to the dynamics, we assume that we can only measure variable $u$
at discrete and sparse points. More precisely, we measure the proximity between two state $(u,v)$ and $(\hat u,\hat v)$ by
\begin{equation}
\int_\Omega |u(x)-\hat u(x)|^2\phi(x) \id x,
\end{equation}
which only relies on the variable $u$, with $\phi$ being the optimal norm density to be learned using the SPML algorithm.

In steps 1 and 2 of Algorithm~\ref{alg:SPML}, we generate the labeled precomputed library with the initial conditions $\left(u_0^{(i)},v_0^{(i)}\right)$. 
To this end, we first generate the random smooth functions,
\begin{align} \label{eq:IC FHN}
    u(x,0)  & = u_2 + \frac{1}{K} \sum_{k=1}^{K} \Big[a_k\cos(k\pi x) \nonumber \\
     &\quad + b_k \sin\left(\frac{2k-1}{2}\pi x\right)\Big],\nonumber\\
    v(x,0) & = 0 
\end{align} 
where $K=22$, $a_k$ and $b_k$ are random numbers drawn from a standard normal distribution, and $u_2$ is the unstable steady state, defined earlier.

\begin{figure}
    \centering  
    \includegraphics[width=0.6\linewidth]{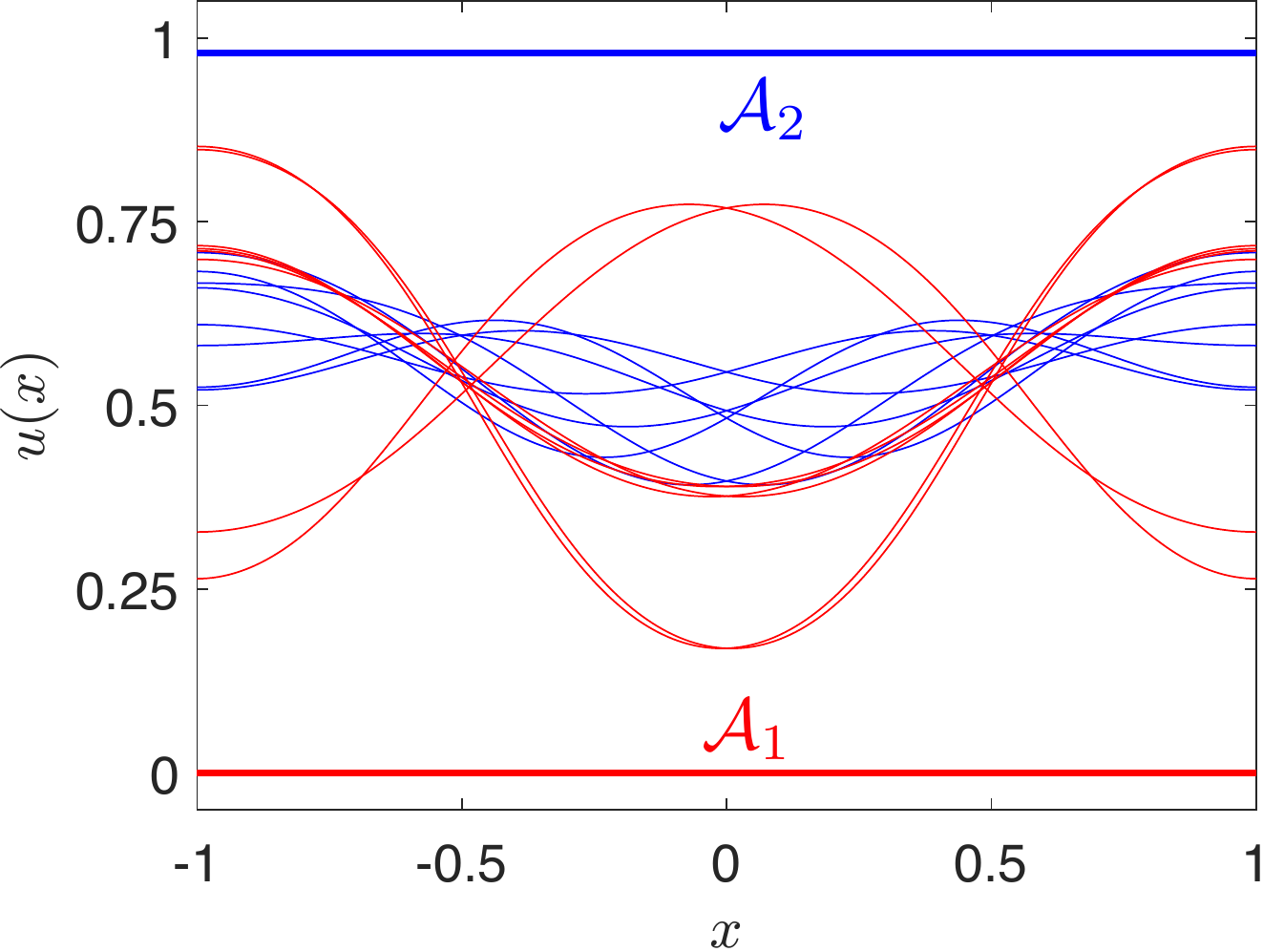}
    \caption{Two steady state attractors $\mathcal A_1$ and $\mathcal A_2$, and sixteen labeled initial conditions $u_0(x)=u(x,t_0)$. The initial conditions are obtained from
    	the random functions~\eqref{eq:IC FHN} after evolving them for $t_0=10$ time units under FitzHugh--Nagumo equation~\eqref{FHN_eqn}. 
    Each state is colored according to the attractor to which it converges asymptotically.}
\label{fig:FHN_initialconditions}
\end{figure}

We evolve the random functions~\eqref{eq:IC FHN} under the FitzHugh--Nagumo equation for a short period of ten time units to obtain the initial conditions $u(x,t_0)=u_0(x)$, where $t_0=10$.
The resulting initial conditions $u_0$ and their attractors are shown in figure~\ref{fig:FHN_initialconditions}.
We use the evolved states $u_0(x)$ as the labeled data. 
This is in line with applications, where state observations are made after the state has evolved for a certain period of time under the system dynamics.

In step 2, we label the initial conditions. 
To this end, we evolve the initial conditions until they converge to their asymptotic state.
If initial condition $u_0^{(i)}$ converges to attractor $\mathcal A_\alpha$, with $\alpha\in\{1,2\}$, we assign the label $\ell(i)=\alpha$ to that initial condition. 
On average it takes about $100$ time units for the initial conditions to converge to an attractor.

We note that, similar to the reaction-diffusion equation, the FitzHugh--Nagumo equation~\eqref{FHN_eqn} is equivariant under the transformation $x\mapsto -x$.
As a result if $u(x,t)$ is a solution to the FitzHugh--Nagumo equation corresponding to the initial condition $u_0(x)$, then 
$u(-x,t)$ is the solution corresponding to the initial condition $u_0(-x)$. In particular, this implies that if $u_0(x)$ converges 
to attractor $\mathcal A_1$ (respectively, $\mathcal A_2$), then the initial condition $u_0(-x)$ also converges to $\mathcal A_1$ (respectively, $\mathcal A_2$).
Using this equivariance, one can double the size of the precomputed library at no additional computational cost. 
The results reported below include these symmetric copies.

\begin{figure}
	\centering
	\subfigure[\label{fig:5labels_FHN} ]{\includegraphics[width=.49\textwidth]{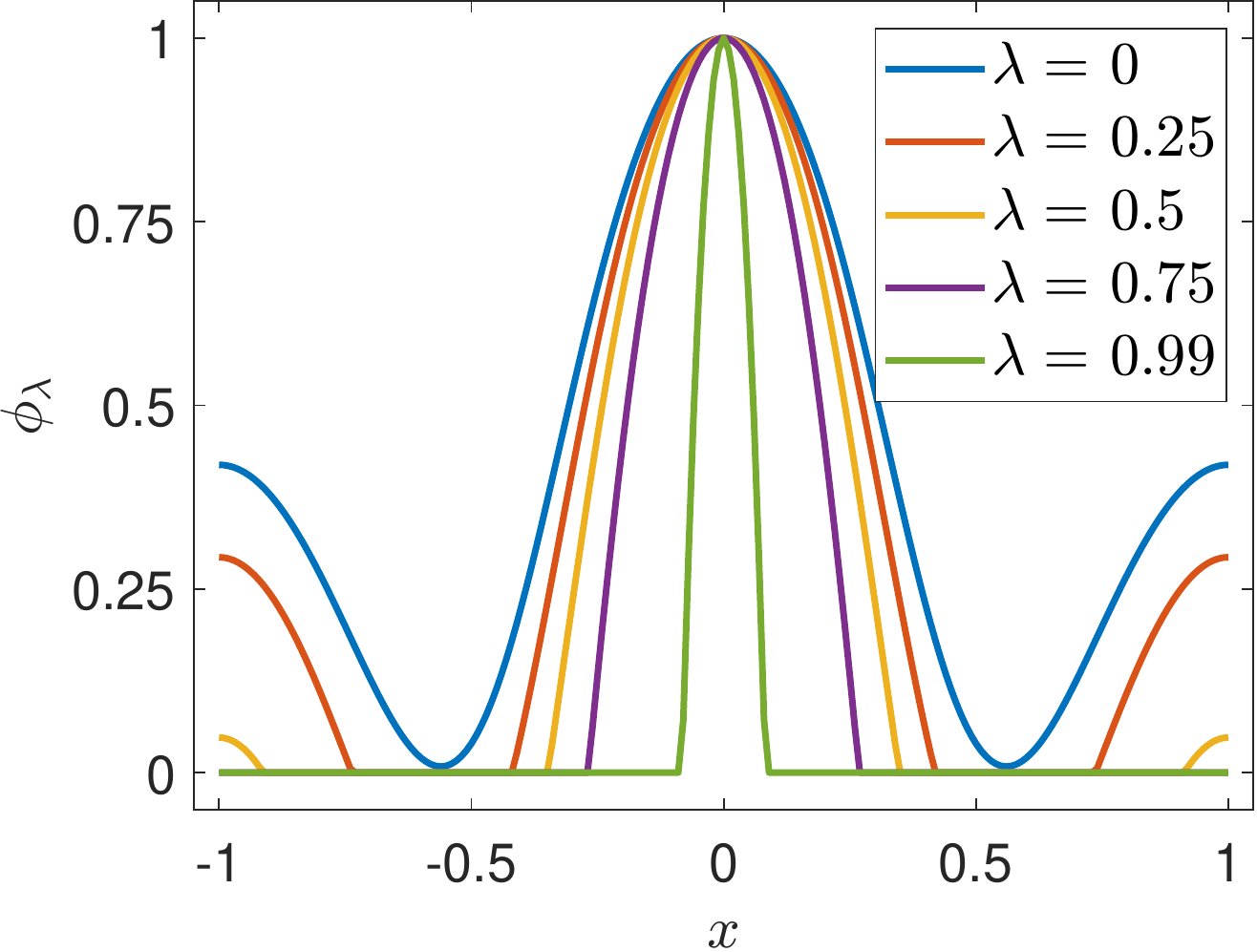}}
	\subfigure[\label{fig:10labels_FHN}]{\includegraphics[width=.49\textwidth]{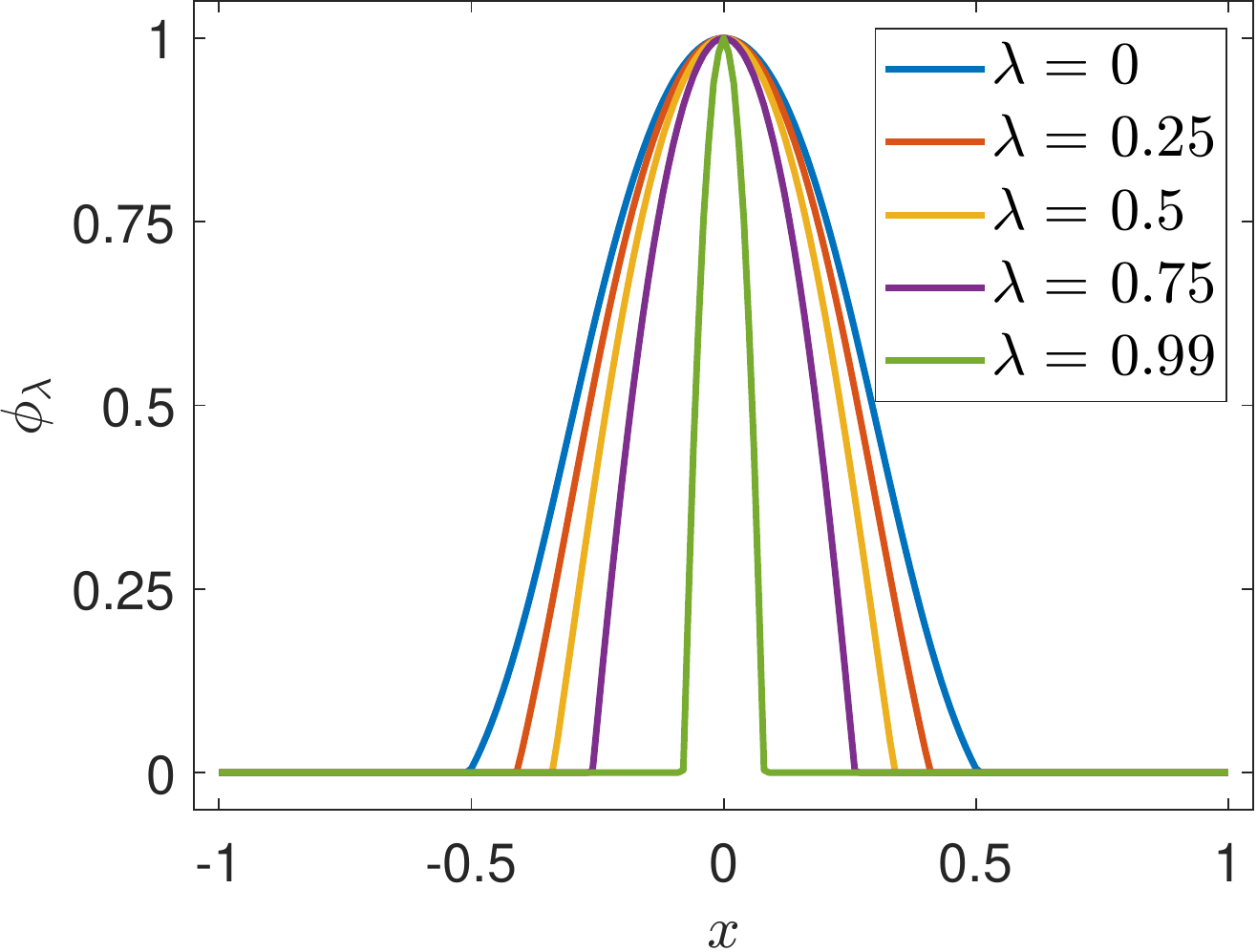}}
	\caption{The learned density $\phi_{\lambda}$ with $\a= 1$ and various values of $\lambda$. 
    The densities are learned from (a) 10 labeled data per attractor. 
    (b) 20 labeled data per attractor. 
The densities are normalized by their maximum value, $\max_{x\in\Omega} \phi_{\lambda}(x)$.} 
	\label{fig:phi_FHN}
\end{figure}

In step 3, we solve the optimization problem~\eqref{original_opt} in the manner discussed in Section~\ref{sec:opt}. 
We choose the regularization parameter $\alpha=1$ and vary the sparsity promoting parameter $\lambda$. 
Figure~\ref{fig:phi_FHN} shows the corresponding optimal densities $\phi_\lambda$ computed from a library of initial conditions $u_0^{(i)}$ consisting of $10$ and $20$ labeled data per attractor.
With no sparsity promotion, i.e. $\lambda=0$, the optimal density $\phi_0$ is non-zero over large portions of the domain $[-1,1]$. 
Therefore, the corresponding norm (cf. equation~\eqref{eq:phi_norm}) cannot be used to quantify the proximity of states when only sparse measurements are available.
We refer to $\phi_0$, the optimal density at $\lambda=0$, as the \emph{non-sparse optimal density} for the FitzHugh--Nagumo equation.

The optimal density $\phi_\lambda$ gradually concentrates around $ x = 0$ as $\lambda$ increases. 
Taking this observation to its logical limit, $\phi_s = \delta(x)$ is  the \emph{sparse optimal density} for the FitzHugh--Nagumo equation. 
The pseudo-metric corresponding to $\phi_s$ is 
\begin{align}\label{eq:snorm_FHN}
    \|u-\hat u\|_{\phi_s} & = \left(u(0)-\hat u(0)\right)^2,
\end{align}
which measures the pseudo-distance between two states $(u,v)$ and $(\hat u,\hat v)$, when only measurements of the 
$u$ variable are possible.
The SPML optimization finds a single location for the optimal sensor location so the pseudo-metric is computable from single point measurements. 
In other words, for the FitzHugh--Nagumo equation, the sparsity-promoting optimization~\eqref{original_opt} informs the optimal sensor placement $\mathcal G_s = \{0\}$ in addition  to returning an optimal norm for quantifying the proximity of system states.
Therefore, in step 4 of Algorithm~\ref{alg:SPML}, we only make state measurements at $x=0$. 

\begin{figure}
    \centering
    \includegraphics[width=0.8\textwidth]{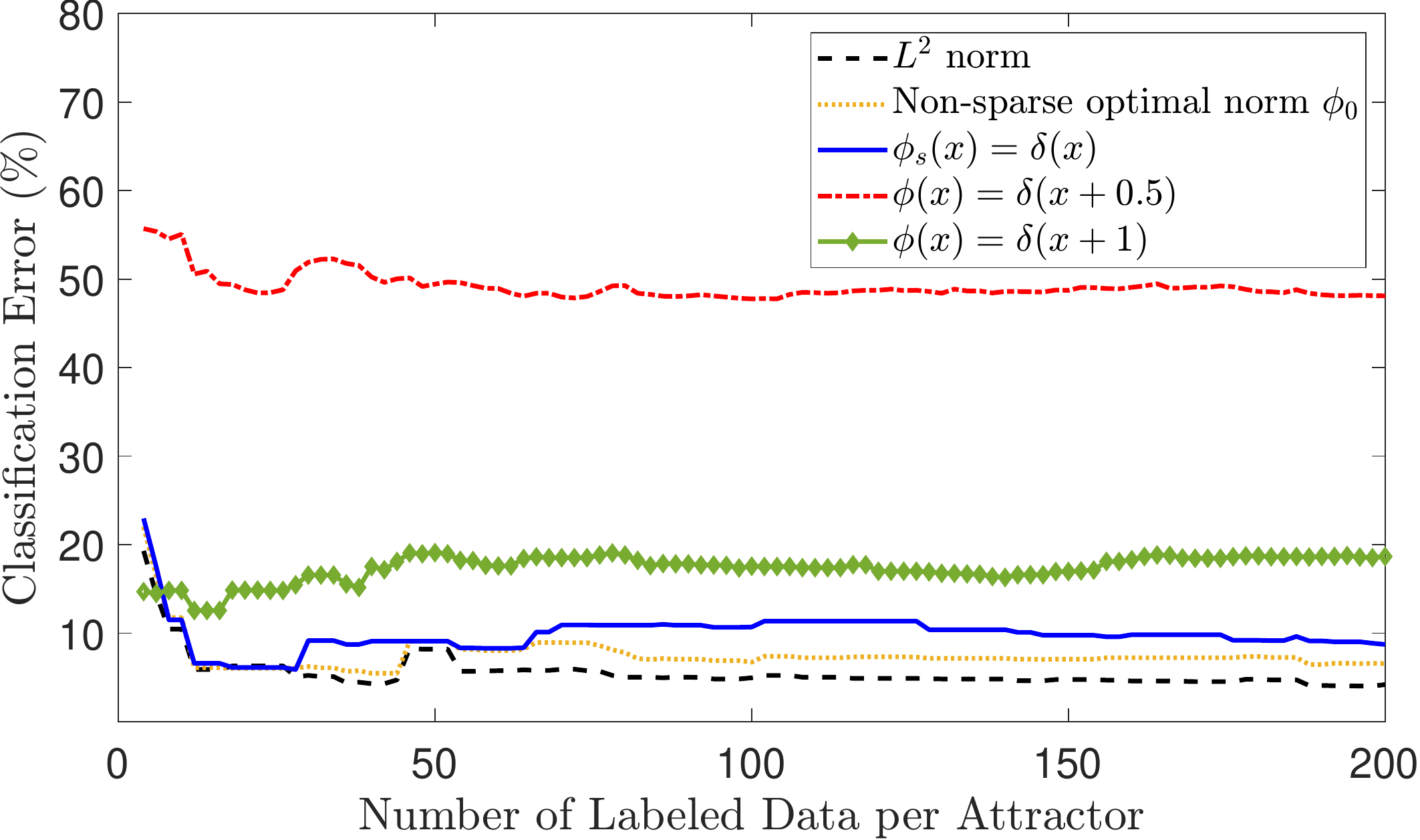}
    \caption{Classification errors as a function of the number of labeled data for the FitzHugh--Nagumo system. Note that the $L^2$ norm and non-sparse optimal norm
    $\|\cdot\|_{\phi_0}$, although performing well, are not computable from sparse measurements.}%
    \label{fig:FHN errors}
\end{figure}

We predict asymptotic behavior of new out-of-library initial conditions $u_0(x)$ uses the nearest neighbor clustering,
where the distances between states are measured using the sparse optimal norm~\eqref{eq:snorm_FHN}.
Figure~\ref{fig:FHN errors} shows the classification errors for various numbers of labeled data. 
The errors were calculated from $3000$ out-of-library initial conditions. 
For comparison, we also show the classification errors from the usual $L^2$ norm,
the non-sparse optimal norm $\|\cdot\|_{\phi_0}$, and
sparse norms computed from single sensors located at the sub-optimal points $x=-0.5$ and $x=-1$.

The classification errors corresponding to the $L^2$ norm, the non-sparse optimal norm, and the optimal sparse optimal norm  exhibit a similar trend:
as the number of labeled data $n_\ell$ increases, the classification errors generally decrease, and saturate around $5-10\%$ error
for a moderate number of labeled data.
In contrast to the reaction-diffusion equation, the $L^2$ norm and the non-sparse optimal norm outperform the sparse optimal norm.
We emphasize, however, that the $L^2$ norm and non-sparse optimal norm are not computable from sparse measurements.

We also compare the optimal sparse norm to sub-optimal sparse norms, computed from on-point measurements.
The sparse norm evaluated from measurements at $x=-0.5$ performs significantly worse than any other norm. 
The sparse norm evaluated from measurements at  $x=-1$, although performing well for small number of labeled data (less than 6 labeled data per attractor),
also loses to the optimal sparse norm for larger number of labeled data.

\section{Conclusions} \label{sec:Conclusion}
We proposed a data-driven method for predicting the asymptotic limit of multistable systems based on sparse 
measurements of their initial condition. This method falls in the framework of semi-supervised 
classification where the labeled data are provided through a library of initial conditions that are evolved numerically
to determine their asymptotic behavior. Given sparse measurements of an out-of-library initial condition, its attractor is determined
based on its proximity to the labeled data.

The choice of the metric for quantifying this proximity is crucial.
We formulated and solved an optimization problem in order to learn an optimal metric from the library of labeled data.
The optimization problem is set up to ensure that the resulting optimal metric (i) is consistent with the labeled data, i.e., the initial conditions that converge to 
the same attractor are relatively close in this metric, and (ii) is computable from sparse measurements.
We introduced two regularizations to the optimization problem. An $L^2$ regularization to rule out degenerate minimizers, and an $L^1$ regularization
to ensure that the resulting metric is computable from sparse measurements.

We demonstrated the application of our method on {\cb two examples: a reaction-diffusion equation with four asymptotically stable steady states and the FitzHugh--Nagumo equation with two asymptotically stable steady states.}
For the reaction-diffusion equation, the learned optimal norm density peaks at two points and is practically zero everywhere else. As such, it informs the optimal sensor placement
for gathering measurements. Predictions based on these two-point measurements are $95\%$ accurate when a moderate number of labeled data (e.g., $50$ per attractor) are used. 
{\cb For the FitzHugh--Nagumo equation, the learned optimal norm density peaks at a single point, informing optimal sensor placement at $x=0$. Predictions based on this single point are 93\%  accurate at 28 labeled data per attractor}. 
Our method shifts the bulk of the computational cost (i.e., creating the labeled data and solving the optimization problem) to offline computations that
are carried out only once. As a result, real-time predictions from sparse measurements are feasible in a matter of seconds.

Future work will focus on several extensions and applications of the present framework. 
Our method can be immediately applied to problems that can be formulated as a semi-supervised classification problem, e.g.,
prediction of extreme events in PDEs~\cite{Farazmand2017,farazmand2019a}.
Although here we presented the results for a one-dimensional PDE, the framework is applicable to higher dimensional problems. 
However, the SPML optimization might become computationally prohibitive in higher dimensions, demanding specialized fast algorithms for its solution.
Finally, here we focused on the class of weighted $L^2$ norms with density $\phi$; however, other types of norms might be more suitable for a specific problem and should be explored.

\section*{Acknowledgments}  
We are grateful to Prof. Robert L. Jerrard (University of Toronto) for his comments on gradient flows, and 
to Prof. Ashley P. Willis (University of Sheffield) for recommending the visualization shown in figure~\ref{fig:projTraj}.

\appendix
\section{Proof of Theorem~\ref{thm:opt}}\label{app:rank1}
In this section, we present the proof of Theorem~\ref{thm:opt}. 

Item 1. Convexity of the optimization follows from the fact that the objective functional
\begin{equation}
J(\phi)=S(\phi)+ \a \big(\lambda\|\phi\|_1+(1-\lambda)\|\phi\|_2  \big), 
\end{equation}
and the constraints are convex. The convexity of the objective functional $J$ follows from linear dependence of the similar sum $S(\phi)$ on $\phi$, and the fact
that $L^p$ norms, with $p\geq 1$, are convex. 

Item 2. Assume that the minimizer satisfies $D(\phi)= D_0$, for some $D_0>1$. Define $\hat \phi = \phi/D_0$ and note that $D(\hat \phi)=1$, $\hat \phi>0$, and 
$J(\hat \phi) = J(\phi)/D_0<J(\phi)$. This contradicts $\phi$ being the minimizer. Therefore, we must have $D(\phi)=1$.

Item 3. Let $\tilde\phi$ denote the minimizer of the optimization problem~\eqref{original_opt}, and $\hat \phi$
denote the minimizer of the scaled optimizaton problem where the constraint~\eqref{second} is replaced by $D(\phi)\geq c$, for some $c>0$. We prove that $\hat \phi = c \tilde \phi$. 

Assume that $\hat \phi\neq c\tilde \phi$. Since $\hat \phi$ is the minimizer of the scaled problem, we have
$$J(\hat \phi)< J(c\tilde\phi),\quad D(\hat \phi)\geq c,\quad \hat \phi>0.$$
Since $J$ is a homogeneous functional of degree one, and $D$ is linear in $\phi$, we obtain
$$J(\hat \phi/c)< J(\tilde\phi),\quad D(\hat \phi/c)\geq 1,\quad \hat \phi/c>0.$$
Therefore, $\hat \phi/c$ satisfies the constraints of~\eqref{original_opt} and the objective functional $J(\hat\phi/c)$
is smaller than $J(\tilde \phi)$, which contradicts $\tilde \phi$ being the minimizer. As a result, we must have
$\hat \phi = c \tilde \phi$.

Item 4. Proof of this property takes considerably more work.
This property shows that the $L^2$ penalization $\|\phi\|_2$ is necessary to rule out degenerate solutions of the optimization problem~\eqref{original_opt}. 
More specifically, without this term, the optimal solution can be non-zero only at a single point. We show this for the discretized version of the 
optimization problem, using collocation points. We discretize the domain $\Omega$ into $N$ disjoint sets $\Omega_k$, such that
$\Omega=\overline{\cup_{k=1}^N\Omega_k}$.
We approximate the integrals with the Riemann sum,
\[\int_\Omega f(x) \, \id x \simeq \sum_{k} f(x_k)\, \mbox{vol} (\Omega_k),\]
where $x_k\in\Omega_k$ and $\mbox{vol} (\Omega_k)$ denotes the volume of the set $\Omega_k$. 

To simplify notation further, we define, 
$$\phi_k :=\phi(x_k)\mbox{vol} (\Omega_k), \quad u^{ij}_k =\left (u^{(i)}_0(x_k)-u^{(i)}_0(x_k)\right)^2,$$
where the weight $\mbox{vol} (\Omega_k)$ is absorbed into the definition of $\phi_k$.
With this discretization, the optimization problem~\eqref{original_opt} with no regularization ($\alpha=0$) becomes,
\begin{subequations}
\begin{equation} 
    \min_{\phi_k}\sum_{\sS} \sum_{k} u_k^{ij}\phi_k,
\end{equation}
\begin{equation}
\sum_{\sD}  \sum_{k} u_k^{ij}\phi_k \geq 1,
\end{equation}
\begin{equation}
\phi_k \geq 0. 
\end{equation}
\end{subequations}

We define the augmented Lagrangian,
	\begin{align}
        J(\phi,\mu, \gamma, \etab,\bxi)  =& \sum_S\sum_k u^{ij}_k \phi_k \nonumber\\  
        & + \mu \left[ \sum_{\mathcal D} \sum_k u^{ij}_k\phi_k-1-\gamma^2\right] \nonumber\\ 
        &+ \sum_k \eta_k (\phi_k-\xi_k^2)         
	\end{align}
	where $\gamma$, $\mu$, $\bxi=(\xi_1,\cdots,\xi_N)$, and $\pmb\eta = (\eta_1,\cdots,\eta_N)$ are Lagrange multipliers.
    The gradient of $J$ with respect to all variables must vanish at the minimizer, implying
	\begin{subequations}
		\begin{equation}
		\partial_{\phi_\ell}J = \sum_\sS u^{ij}_\ell+\eta_\ell + {\mu} \sum_D {u^{ij}_\ell} =0, 
		\label{eq:dJdA}
		\end{equation}
		\begin{equation}
		\partial_\mu J = \sum_\sD \sum_k u^{ij}_k\phi_k-1-\gamma^2=0,
		\end{equation}
		\begin{equation}
		\partial_\gamma J = -2\mu \gamma =0,
		\end{equation}
		\begin{equation}
		\partial_{\eta_\ell}J= \phi_\ell - \xi_\ell^2 =0,
		\end{equation}
		\begin{equation}
		\partial_{\xi_\ell} J= -2\eta_\ell\xi_\ell = 0,
		\end{equation}
	\end{subequations}
for all $\ell=1,2,\ldots, N$.
These set of equations, imply the following:
\begin{enumerate}
	\item We must have $\gamma=0$. Otherwise, we would have $\mu=0$ which implies $\eta_\ell =-\sum_{\mathcal S} u^{ij}_\ell\neq 0$ for all $\ell$. This implies $\xi_\ell =0$ for all $\ell$, resulting in $\phi\equiv 0$.
	\item As a consequence of $\gamma=0$, the minimizer must satisfy $\sum_{\mathcal D}\sum_k u^{ij}_k\phi_k =1$.
	\item If $\phi_\ell\neq 0$, then $\eta_\ell =0$. This follows from $\phi_\ell = \xi_\ell^2$ and $\eta_\ell\xi_\ell=0$. 
\end{enumerate}
The last observation implies that the minimizer $\phi$ is degenerate, in the sense that generally $\phi_\ell$ can be nonzero exactly at one point. Proof of this is by contradiction. Assume that $\phi_\ell$ is nonzero for two indices $\ell\in\{1,\ldots, N\}$. 
Without loss of generality, assume that these indices are $\ell=1,2$ so that $\phi_1\neq 0$ and $\phi_2\neq 0$. As a result, $\eta_1=\eta_2=0$.
This together with equation~\eqref{eq:dJdA} implies that
\begin{equation}
    \mu = - \frac{\sum_\sS u^{ij}_1}{\sum_\sD u^{ij}_1}=- \frac{\sum_\sS u^{ij}_2}{\sum_\sD u^{ij}_2}
\end{equation}
This is a contradiction since the last two expressions are not equal for generic labeled data $u_0^{(i)}$.
The proof for $\lambda=1$ is similar and therefore omitted for brevity.

Xing et al.~\cite{Xing2003} use the slightly different constraint $\hat D(\phi)\geq 1$, where 
\begin{equation}
\hat D (\phi) = \sum_\sD \|u_0^{(i)}-u_0^{(j)}\|_\phi.
\end{equation}
Note that this is different from the dissimilar sum~\eqref{eq:dissum} we use, in that the norms in the sum are not squared.
Xing et al.~\cite{Xing2003} imply that the constraint $\hat D(\phi)\geq 1$ rules out the degenerate rank-1 minimizers. 
This statement, however, is not correct. 
A similar analysis, as the one presented above, shows that if $\phi_\ell$ is nonzero for two indices, say $\ell =1,2$, then
we must have
\begin{align}
    \mu  & = -2 \frac{\sum_\sS u^{ij}_1  }{\sum_\sD \frac{u^{ij}_1}{(\sum_{k=1}^N u^{ij}_k\phi_k)^{1/2}}} \nonumber\\ & = -2 \frac{\sum_\sS u^{ij}_2}{\sum_\sD \frac{u^{ij}_2}{(\sum_{k=1}^N u^{ij}_k\phi_k)^{1/2}}}.
\end{align}
Again this is a contradiction because $\mu$ assumes two values that are generically not equal.

Now we show that the $L^2$ regularization allows for non-degenerate minimizers. 
Consider the case where $\a >0, \lambda=0$, 
which corresponds to the optimization problem~\eqref{original_opt} with only $L^2$ regularization. The case $0<\lambda<1$ which includes both 
$L^1$ and $L^2$ regularizations is similar.
Without loss of generality we assume $\a=1$. The augmented Lagrangian then reads
	\begin{align}
        J(\phi,\mu, \gamma, \etab,\bxi) =&  \sum_S\sum_k u^{ij}_k \phi_k  
         +\sum_k \eta_k (\phi_k-\xi_k^2) \nonumber\\& +\mu \left[ \sum_D \sum_k u^{ij}_k\phi_k-1-\gamma^2\right]\nonumber\\  
        & +     \left(\sum_{k}\phi_k^2\right)^{1/2}
	\end{align} 
    The gradient of $J$ with respect to all variables must vanish at the minimizer.
	Therefore, 
	\begin{subequations}
		\begin{equation}\label{eqn:better_J}
            \partial_{\phi_\ell}J = \sum_\sS u^{ij}_\ell + {\mu} \sum_{\sD} {u^{ij}_\ell}+\eta_\ell + \frac{\phi_\ell}{\|\phi\|_2} = 0, 
		\end{equation}
		\begin{equation}
		\partial_\mu J = \sum_\sD \sum_k u^{ij}_k\phi_k-1-\gamma^2=0,
		\end{equation}
		\begin{equation}
		\partial_\gamma J = -2\mu \gamma =0,
		\end{equation}
		\begin{equation}
		\partial_{\eta_\ell}J= \phi_\ell - \xi_\ell^2 =0,
		\end{equation}
		\begin{equation}
		\partial_{\xi_\ell} J= -2\eta_\ell\xi_\ell = 0,
		\end{equation}   
	\end{subequations}
for all $\ell=1,2,\ldots, N$.
The addition of the $\frac{\phi_\ell}{\|\phi\|_2}$ allows extra degrees of freedom in the solution to equation~\eqref{eqn:better_J}, so that the solution is not necessarily rank-1.
It is still true that if $\phi_\ell\neq 0$, then we must have $\eta_\ell=0$. But now the $L^2$ regularization and equation~\eqref{eqn:better_J} imply
\begin{equation}
\phi_\ell = -\|\phi\|_2\left( \sum_\sS u^{ij}_\ell + {\mu} \sum_{\sD} {u^{ij}_\ell}\right),
\end{equation}
for non-zero elements of $\phi$ which does not lead to a contradiction.

\section{Gradient flow structure}
\label{app:inorm}
In this section, we analyze the gradient flow structure of the reaction-diffusion equation~\eqref{generic_RDE} and show that it has intrinsic inner product structure.
Consider the energy functional,
\begin{equation}
    I[u] = \int_\Omega \left[ \frac{\nu}{2}|\grad  u|^2 - F(u) \right] w \, \id x,
\end{equation}
and weighted inner product,
\begin{equation}
    \left<u,v \right>_w = \int_\Omega u(x)\cdot v(x)w(x)\, \id x.
\end{equation}
The resulting gradient flow is given by $\frac{\d u}{\d t} = -\grad_H I(u)$, where the gradient $\grad_H I$ is a map that satisfies $\frac{\id}{\id \varepsilon} I[u+\varepsilon v] \Big|_{\varepsilon=0} = \left<\grad_H I(u), v \right>_w$, for all $v\in H$.
Using integration by parts, the variational derivative can be written as,  
   \begin{align*}
       \frac{\id}{\id \varepsilon}\Bigg|_{\varepsilon=0} & I[u+\varepsilon v] \\
       & = \frac{\id}{\id \varepsilon}\Bigg|_{\varepsilon=0} \int_\Omega \Big[ \frac{\nu}{2}|\grad u +\varepsilon\grad v|^2 \\ 
        &- F(u+\varepsilon v) \Big] w \, \id x \\  
       =& \int_\Omega \left(\nu \grad u\cdot \grad v - F'(u)v\right) w \, \id x \\
       =& \int_{\partial \Omega} \nu vw \grad u\cdot \hat n\, \id\Gamma  \\
        &+\int_\Omega (-\nu \, \operatorname{div}(w \grad u) v -  F'(u)v w )\, \id x\\
       =& \int_\Omega \left( \frac{-\nu}{w} \text{div}(w \grad u) - F'(u) \right) v w \, \id x \\
        & = \left< -\frac{\nu}{w} \text{div}(w \grad u) - F'(u) ,v \right>_w.
        \end{align*}
Note that the integral over the boundary $\partial \Omega$ vanishes because of the Neumann boundary condition $\hat n\cdot \nabla u=0$, 
where $\hat n$ is the unit normal to the boundary.

Therefore, we have 
\begin{equation}
\nabla _H I(u) = -\frac{\nu}{w} \nabla\cdot(w \grad u) - F'(u),
\end{equation}
and the gradient flow coincides with the reaction-diffusion equation,
\begin{equation*}
    \frac{\d u}{\d t} =  \frac{\nu}{w} \nabla\cdot (w \nabla u)+F'(u).
\end{equation*}
The inner product $\langle\cdot,\cdot\rangle_w$ is intrinsic to this PDE since choosing a different inner product alters the gradient $\nabla _H I(u)$, resulting
in a PDE that is different from the reaction-diffusion equation~\eqref{generic_RDE}.



\end{document}